\documentclass{article}
\usepackage{authblk}
\usepackage[T1]{fontenc}
\usepackage[utf8]{inputenc}
\usepackage{lmodern} 
\usepackage[english]{babel}
\usepackage{hyperref}
\hypersetup{
	pdftitle={Shape optimization in the space of piecewise-smooth shapes for the Bingham flow variational inequality},
	pdfauthor={Tim Suchan, Volker Schulz and Kathrin Welker}
}

\usepackage{amsmath,amssymb,amscd,latexsym,amsthm}
\usepackage{upgreek}
\usepackage{amsfonts}

\usepackage{mathtools}
\usepackage{bm}

\newtheorem{remark}{Remark}

\newenvironment{acknowledgements}
{\par\footnotesize}
{\par\addvspace{\bigskipamount}}

\usepackage{cleveref}

\usepackage{algorithm}
\usepackage[noend]{algpseudocode}

\usepackage{color}
\usepackage{graphicx}
\graphicspath{{figures/},{svg-inkscape/},{tikzs/}}
\usepackage{subcaption}
\newlength\figureheight 
\newlength\figurewidth

\newcommand{\includetikz}[1]{\includegraphics{tikzs/#1_compiled}}

\newcommand{\tunderbrace}[2]{\underbrace{#1}_{\textstyle#2}}
\newcommand{\shape}{u}
\newcommand{\admissibleShapeSpace}{M_1^{\mathrm{ad}}}
\newcommand{\compDomain}{\varOmega}
\newcommand{\dx}{\,\mathrm{d} \bm{x}}
\newcommand{\velocity}{\bm{y}}
\newcommand{\pressure}{p}
\newcommand{\adjvelocity}{\tilde{\bm{y}}}
\newcommand{\adjpressure}{\tilde{p}}
\newcommand{\perturbvelocity}{\hat{\bm{y}}}
\newcommand{\perturbpressure}{\hat{p}}
\newcommand{\perturbvelocitytwo}{\adjvelocity}%
\newcommand{\perturbpressuretwo}{\adjpressure}%
\newcommand{\viscosity}{\mu}
\newcommand{\density}{\rho}
\newcommand{\yieldStress}{g}
\newcommand{\regParam}{\gamma}
\newcommand{\regParamMax}{\delta}
\newcommand{\regNonDifferentiability}{\bm{h}_\regParam}
\DeclareMathOperator{\Div}{div}
\newcommand{\stress}{\bm{\sigma}}
\newcommand{\strain}{\bm{\varepsilon}}
\DeclareMathOperator{\vol}{vol}
\DeclareMathOperator{\bary}{bary}
\DeclareMathOperator{\peri}{peri}
\newcommand{\const}{\textrm{const.}}
\newcommand{\meshDeformation}{\bm{V}}
\newcommand{\perturbshape}{\meshDeformation}

\begin{document}

\title{Shape optimization in the space of piecewise-smooth shapes for the Bingham flow variational inequality}

\author{Tim Suchan}
\affil{Helmut-Schmidt-Universität/Universität der Bundeswehr Hamburg,  Holstenhofweg~85, 22043 Hamburg, Germany, \href{mailto:suchan@hsu-hh.de}{\ttfamily suchan@hsu-hh.de}\vspace*{6pt}}

\author{Volker Schulz}
\affil{Universität Trier, Universitätsring 15, 54296 Trier, Germany, \href{mailto:volker.schulz@uni-trier.de}{\ttfamily volker.schulz@uni-trier.de}\vspace*{6pt}}

\author{Kathrin Welker}
\affil{TU Bergakademie Freiberg,  Akademiestr.~6, 09599 Freiberg, Germany, \href{mailto:Kathrin.Welker@math.tu-freiberg.de}{\ttfamily Kathrin.Welker@math.tu-freiberg.de}}

\date{}

\maketitle

\begin{abstract}
This paper sets up an approach for shape optimization problems constrained by variational inequalities (VI) in an appropriate shape space. In contrast to classical VI, where no explicit dependence on the domain is given, VI constrained shape optimization problems are in particular highly challenging because of two main reasons: Firstly, one needs to operate in inherently non-linear, non-convex and infinite-dimensional shape spaces. Secondly, the problem cannot be solved directly without any regularization techniques in general because, e.g., one cannot expect the existence of the shape derivative for an arbitrary shape functional depending on solutions to VI. This paper introduces a specific shape manifold and presents an optimization technique to handle the non-differentiabilities on this shape manifold. In particular, we formulate an optimization system based on Gâteaux semiderivatives and Eulerian derivatives for a shape optimization problem constrained by the Bingham flow variational inequality. Numerical results show the applicability and efficiency of the proposed approach.
\end{abstract}

\section{Introduction}
\label{sec:Intro}

We focus on an objective functional that is supposed to be minimized with respect to a shape.
Finding a correct model to describe the set of shapes is one of the main challenges in shape optimization.
Since algorithmic ideas from \cite{Absil} can be combined with approaches from differential geometry and a Riemannian metric yields the distances of two shapes, it is attractive to optimize on Riemannian manifolds from a theoretical and computational point of view (cf., e.g., \cite{LoayzaGeiersbachWelkerHandbook}).
If one concentrates on one-dimensional shapes, the shape space
$B_e(S^1, \mathbb{R}^2)$
briefly investigated in~\cite{MichorMumford1} is an important example of a smooth manifold allowing a Riemannian structure, where the term smooth shall refer to infinite differentiability in this paper. This shape space is considered in recent publications (cf., e.g., \cite{Geiersbach2021,geiersbach2023stochastic,Schulz,Schulz2016a}),  but  $B_e(S^1, \mathbb{R}^2)$ is in general not sufficient to carry out optimization algorithms on non-smooth shapes.
Recently, the optimization of piecewise-smooth shapes has been addressed  \cite{GeiersbachSuchanWelker_NS,Pryymak2023}. 
In \cite{Pryymak2023}, a shape is seen as a point on an abstract manifold so that a collection of shapes can be viewed as a vector of points on a product manifold.
More precisely, a space containing shapes in $\mathbb{R}^2$ that can be identified with a Riemannian product manifold but at the same time admits piecewise-smooth curves as elements is constructed.
In this paper, we focus on this space of piecewise-smooth shapes. This novel shape space was  used in an optimization procedure in \cite{GeiersbachSuchanWelker_NS}, in which also a framework for the handling of piecewise-smooth shapes is provided.
Since this shape space is a product manifold, optimization techniques on it base on the optimization of multiple shapes firstly presented in~\cite{LoayzaGeiersbachWelkerHandbook} and applied in \cite{geiersbach2023stochastic,Pryymak2023}.

It is common to encounter optimization problems with a partial differential equation (PDE) constraint that describes the physical behavior of the system. These PDE depend on the field of application, e.g., structural mechanics \cite{Allaire2004,Haslinger2003,Upadhyay2021}, acoustics \cite{Kapellos2019,Schmidt2016} and fluid mechanics \cite{Bletsos2023,Haslinger2003,Kuehl2022} and can even include multiple of these fields \cite{Haubner2020,Wiegard2021}. However, the physical behavior cannot always be described by a PDE. Obscacle-like behavior like a maximal value for the unknown of the constraint, as encountered in contact problems, can be described by variational inequalities (VI) of the first kind \cite{Duvaut1976,Haslinger2003}, and a switch in behavior depending on the unknown of the constraint, e.g. friction problems or a change in fluid behavior are modelled by VI of the second kind \cite{Duvaut1976,Fuchs2000}.
In this paper, we concentrate on a VI of second kind.

There is a sizeable literature on shape optimization constrained by partial differential equations. In contrast to that, a much smaller number publications consider VI as the state model of a shape optimization problem. In \cite[Chap.~4]{SokoZol}, shape derivatives of elliptic VI problems are presented. Also, \cite{LR-1991a} presents existence results for shape optimization problems which can be reformulated as optimal control problems, whereas \cite{DM-1998,G-2001} show existence of solutions in a more general set-up. In \cite{Myslinski-2007}, level-set methods are proposed and applied to graph-like two-dimensional problems. Moreover, \cite{HL-2011} presents a regularization approach to the computation of shape and topological derivatives in the context of elliptic VI and, thus, circumvents the numerical problems in \cite[Chap.~4]{SokoZol}. In \cite{Sturm-2016} shape optimization for VI in the context of damage models are treated. Directional derivatives for shape optimization problems with VI are discussed in \cite{Kovtunenko2023-qg}. In \cite{Luft2020,SchulzWelker-2021}, necessary conditions of adjoint type are derived for shape optimization in the context of the obstacle problem, as the limit of a regularization approach.  

In this paper, we study Bingham flow as a variational inequality constraint, which models a sudden change from solid to fluid behavior depending on if the stresses in the fluid are below or above a threshold, to a shape optimization problem and derive necessary conditions on the basis of Gâteaux semiderivatives. Approximations to these necessary conditions are employed in an optimization algorithm and numerical results are presented  for the first time---to the knowledge of the authors.
The paper is structured as follows. In \Cref{sec:ModelFormulation}, we introduce the setting for optimization on Riemannian manifolds and introduce a specific manifold (space of piecewise-smooth shapes) considered in this paper and  the model problem of shape optimization for Bingham fluids. In \Cref{sec:ApplicationSemiderivativeApproach} we formulate the optimality system and present an optimization technique to handle the non-differentiabilities.  \Cref{sec:Num} then describes the numerical implementation of the optimization problem and shows numerical results.

\section{Model formulation}
\label{sec:ModelFormulation}

In \Cref{sub:ShapeSpace}, we introduce a space containing shapes in $\mathbb{R}^2$ that
admits piece\-wise-smooth curves as elements and can be identified with a Riemannian product manifold. This shape space was first defined in \cite{Pryymak2023} and considered in an optimization procedure in \cite{GeiersbachSuchanWelker_NS}.
In \Cref{sub:Model}, we then formulate the shape optimization model that we concentrate on in this paper.

\subsection{Space of piecewise-smooth shapes}
\label{sub:ShapeSpace}

We consider manifolds $\mathcal{U}_1, \dots, \mathcal{U}_N$, each equipped with a Riemannian metrics $G^{i} = (G_u^{i})_{u \in \mathcal{U}_{i}}$, $i=1,\dots,N$. 
The Riemannian product manifold $\mathcal{U}^N$ is then defined by
$  \mathcal{U}^N \coloneqq \mathcal{U}_1 \times \cdots \times \mathcal{U}_N= \prod_{i=1}^{N}\mathcal{U}_{i}$.
In the following, for a point $u_i\in \mathcal{U}_i$, the tangent space of $\mathcal{U}_i$ at $u_i$ is denoted by $T_{u_i} \mathcal{U}_i$, $i=1,\dots,N$. 
As described in  \cite{LoayzaGeiersbachWelkerHandbook}, the tangent space
$T_{\bm{u}}\mathcal{U}^N$ can be identified with the product of tangent spaces $T_{u_1}\mathcal{U}_1\times\dots\times T_{u_N}\mathcal{U}_N$ via 
$
T_{\bm{u}} \mathcal{U}^N \cong T_{u_1} \mathcal{U}_1 \times \dots \times T_{u_N} \mathcal{U}_N.
$
The product metric $\mathcal{G}^N$ to $\mathcal{U}^N$ can be defined by
$\mathcal{G}^N=(\mathcal{G}^N_{\bm{u}})_{\bm{u}\in\mathcal{U}^N}$ with
\begin{equation*}
\mathcal{G}^N_{\bm{u}}(\bm{v},\bm{w}) =\sum_{i=1}^{N} \mathcal{G}_{\pi_i(u)}^{i}(\pi_{i_\ast}\bm{v},\pi_{i_\ast}\bm{w})\qquad\forall \,  
\bm{v},\bm{w} \in T_u \mathcal{U}^N,
\end{equation*}
and $\pi_{i}\colon \mathcal{U}^N\to \mathcal{U}_{i}$, $i=1, \dots, N$, correspond to canonical projections and $\pi_{i_\ast}$ denotes the pushforward\footnote{The derivative of a mapping $f\colon M\rightarrow S$ between two differentiable manifolds $M$ and $S$ is defined using the pushforward.  In a point $u\in M$, it is defined by $(f_*)_u \colon T_u M \rightarrow T_{f(u)} S$ with $(f_*)_u(c):= \frac{\textup{d}}{\textup{d}t} f(c(t))|_{t=0} = (f \circ c)'(0),$ where $c\colon I \to M$, $I \subset \mathbb{R}$ is a differentiable curve with $c(0)=u$ and $c'(0) \in T_u M$.} of $\pi_{i'}$
(cf. \cite{LoayzaGeiersbachWelkerHandbook}). 
In \cite{Pryymak2023},  the $s$-dimensional shape space on $\mathcal{U}^N$ is given by 
\begin{align*}
M_s(\mathcal{U}^N)\coloneqq \Big\lbrace \bm{u}=(u_1, \dots ,u_s) \mid \,& u_i \in  \prod\limits_{l=k_i}^{k_i+n_i-1} \mathcal{U}_{l} \,\, \forall i=1,\ldots,s,\, \sum\limits_{i=1}^s n_i = N \text{ and }\\
&k_1=1, \, k_{i+1}=k_i+n_i
\,\forall i =1,\dots,s-1  \Big\rbrace.
\end{align*}
An element in $M_s(\mathcal{U}^N)$ is defined as a vector of $s$ shapes $u_1,\dots,u_s$, where each shape $u_i$ is an element of the product of $n_i$ smooth manifolds.

\begin{remark}
	As already mentioned in \Cref{sec:Intro}, the shape space
	$B_e(S^1, \mathbb{R}^2)$ is considered in recent publications but this shape space contains only simple closed smooth curves in $\mathbb{R}^2$. In particular, by definition, it does not contain curves which have kinks. The product shape manifold described above from \cite{Pryymak2023} allows for a larger set of possible shapes and also includes shapes from $B_e(S^1, \mathbb{R}^2)$. In computations, the shape needs to be discretized. Then,  the number $n_i$ of glued-together smooth shapes, i.e., the maximum number of kinks of each single closed shape $u_i$ with $i=1,\ldots,s$, can be chosen to be, e.g., the number of nodes belonging to the discretization of this shape (cf. \cite{GeiersbachSuchanWelker_NS}).  %
\end{remark}

In our application, we focus only on the optimization of one shape, i.e., we set $s=1$ and consider $M_1(\mathcal{U}^N)$.
Moreover, we are interested in optimizing piecewise-smooth single closed shapes, i.e., we restrict the choice of shapes in $M_1(\mathcal{U}^N)$ to  piecewise-smooth shapes that are glued together. More precisely, we assume that a shape $u\in M_1(\mathcal{U}^N)$ is single closed:  $u =(u_1,\dots,u_N) \in   (B_e([0,1], \mathbb{R}^2))^{N}$  with the additional conditions 
\begin{equation}
\label{eq:glue-condition}
\begin{split}
& u\colon [0,1)\to \mathbb{R}^2 \text{ injective with }\\
& u_{h}(1)=u_{h+1}(0) \, \forall h = 1,\dots,N-1 \text{ and } u_{1}(0)=u_{N}(1),
\end{split}
\end{equation}
where 
\begin{align*}
B_e([0,1], \mathbb{R}^2) &\coloneqq\text{Emb}([0,1], \mathbb{R}^2) /\text{Diff}([0,1]),
\end{align*}
with $\text{Emb}([0,1], \mathbb{R}^2)$ denoting the set of all embeddings from $[0,1]$ into $\mathbb{R}^2$ and $\text{Diff}([0,1])$ being the set of all diffeomorphisms of $[0,1]$.
The space $B_e([0,1], \mathbb{R}^2)$ contains all simple open smooth curves in $\mathbb{R}^2$ (cf., e.g., \cite[Chapter 2, Definition 2.22]{Kuehnel}) and is a manifold allowing Riemannian structures (cf. \cite{Michor1980}).
The Riemannian manifold structure can be used to define distances between shapes, simulate how one shape arises from another one, state convergence properties of developed algorithms, etc. In this paper, we concentrate on the Steklov-Poincaré metric due to its computational advantages outlined, e.g., in \cite{Schulz2016,Schulz2016a}.
In the following, we identify a closed curve with kinks with a glued-together curve of open smooth curves. Thus, it is convenient to define the admissible space
\begin{align}
\label{AdShapeSpace}
\admissibleShapeSpace \coloneqq \left\lbrace u \in M_1\!\left(B_e([0,1], \mathbb{R}^2)^N\right) \mid u \text{ with } \eqref{eq:glue-condition} \right\rbrace.
\end{align}

\subsection{The Bingham shape optimization problem}
\label{sub:Model}

A common optimization problem for fluid-mechanical models is the minimization of viscous energy dissipation of stationary fluids, expressed as
\begin{align*}
	\min_{\shape \in \admissibleShapeSpace} \int_{\compDomain} \frac{\viscosity}{2} \nabla \velocity : \nabla \velocity \dx,
\end{align*}
cf., e.g.,~\cite{Mohammadi2009,Schulz2016}. For this publication, we denote the domain in which the fluid is located as $\compDomain \subset \mathbb{R}^n$, $n=2, 3$. The fluid velocity~$\velocity$ together with the pressure~$p$ are usually determined from the solution of a partial differential equation (PDE) in $\compDomain$ like the Navier-Stokes or Stokes equations. For this publication, we are interested in a Bingham fluid behavior as described in \cite{Bingham1922} instead. We use the formulation in \cite{Reyes2010} without and \cite{Reyes2011} with a convective term, which reads
\begin{subequations}
\begin{align}
	- \Div{(\stress(\velocity))} + \rho (\velocity \cdot \nabla) \velocity + \nabla \pressure = \bm{f} &\quad \text{in } \compDomain \\
	\Div{(\velocity)} = 0 &\quad \text{in } \compDomain \\
	\velocity = \velocity^D &\quad \text{on } \Gamma^D \subset \partial \compDomain \\
	\stress(\velocity) \, \bm{n} - \pressure \, \bm{n} = \bm{0} &\quad \text{on } \Gamma^N = \partial \compDomain \setminus \Gamma^D \neq \emptyset \\
	\stress(\velocity) = 2 \viscosity \strain(\velocity) + \yieldStress \frac{\strain(\velocity)}{\left\| \strain(\velocity) \right\|} &\quad \text{if } \strain(\velocity) \neq \bm{0} \\
	\left\| \stress(\velocity) \right\| \leq \yieldStress &\quad \text{if } \strain(\velocity) = \bm{0}.\label{2f}
\end{align}
\end{subequations}
Here, $\mu>0$ denotes the fluid viscosity and $\rho>0$ the fluid density, $g>0$ is the plasticity threshold, $\bm{f} \in H^1(\compDomain, \mathbb{R}^n)$ represents the external volume forces, $\| \cdot \|$ is the Frobenius norm, $\stress(\cdot)$ describes the Cauchy stress tensor and $\strain(\cdot)$ stands for the strain tensor, for which we use the symmetric gradient
\begin{align}
	\label{eqn:symmetricGradient}
	\strain(\velocity) = \frac{1}{2} \left(\nabla \velocity + {\nabla \velocity}^\top \right).
\end{align}

The inequality \eqref{2f} does not define $\stress(\velocity)$ uniquely in the case $\strain(\velocity) = \bm{0}$. This is discussed in detail in \cite{Reyes2011, Reyes2010}, where a regularization approach for the treatment of this nonuniqueness is proposed. Here, we follow this approach and thus use a regularizing functional for some $\gamma \gg 0$
\begin{align}
	\regNonDifferentiability(\strain(\velocity)) = \frac{\yieldStress \regParam \, \strain(\velocity)}{\max(\yieldStress, \regParam\! \left\| \strain(\velocity) \right\|)} = \begin{cases}
		\regParam \strain(\velocity) & \text{if } \yieldStress > \regParam \! \left\| \strain(\velocity) \right\|, \\
		\yieldStress \frac{\strain(\velocity)}{\left\| \strain(\velocity) \right\|} & \text{if } \yieldStress \leq \regParam \! \left\| \strain(\velocity) \right\|.
	\end{cases}
	\label{eqn:regularizedNonDifferentiability}
\end{align}
This enables a reformulation of the Bingham flow problem in the following regularized way by defining $\stress(\velocity) = 2 \viscosity \strain(\velocity) + \regNonDifferentiability(\strain(\velocity))$. We define two Hilbert spaces that include Dirichlet boundary information
\begin{align*}
	&\mathcal{H} \coloneqq \{ \velocity \in H^1(\compDomain, \mathbb{R}^n) \colon \velocity = \velocity^D \text{ on } \Gamma^D \} \text{ and } \\ 
	&\mathcal{H}_0 \coloneqq \{ \velocity \in H^1(\compDomain, \mathbb{R}^n) \colon \velocity = \bm{0} \text{ on } \Gamma^D \}.
\end{align*}
The PDE-constrained optimization problem then reads
\begin{align}
	\label{eqn:OptimizationProblemObjectiveFunctional}
	\min_{\shape \in \admissibleShapeSpace} \tunderbrace{\int_{\compDomain} \frac{\viscosity}{2} \nabla \velocity : \nabla \velocity \dx}{J(\shape,\velocity)}
\end{align}
with the PDE constraint in weak form: Find $\velocity, \pressure \in \mathcal{H} \times L^2(\compDomain)$ s.t.
	\begin{align}
		\label{eqn:OptimizationProblemPDEConstraint}
		\begin{aligned}
		 0 = \int_{\compDomain} &2 \viscosity \, \strain(\velocity) : \strain(\perturbvelocitytwo) + \density ((\velocity \cdot \nabla) \velocity) \cdot \perturbvelocitytwo - \pressure \Div(\perturbvelocitytwo) + \perturbpressuretwo \Div(\velocity) \\
		&- \bm{f} \cdot \perturbvelocitytwo + \regNonDifferentiability(\strain(\velocity)) : \strain(\perturbvelocitytwo) \dx \quad \forall \perturbvelocitytwo, \perturbpressuretwo \in \mathcal{H}_0 \times L^2(\compDomain).
		\end{aligned}
	\end{align}
The residual of the weak form of the PDE constraint then is
\begin{align}
	\begin{aligned}
		R(\shape,(\velocity,\pressure),(\perturbvelocitytwo,\perturbpressuretwo)) \coloneqq \int_{\compDomain} &2 \viscosity \, \strain(\velocity) : \strain(\perturbvelocitytwo) + \density ((\velocity \cdot \nabla) \velocity) \cdot \perturbvelocitytwo - \pressure \Div(\perturbvelocitytwo) \\
		& + \perturbpressuretwo \Div(\velocity) - \bm{f} \cdot \perturbvelocitytwo + \regNonDifferentiability(\strain(\velocity)) : \strain(\perturbvelocitytwo) \dx.%
		\label{eqn:OptimizationProblemPDEConstraintResidual}
	\end{aligned}
\end{align}
One method to avoid the occurrence of trivial solutions of the minimization problem \eqref{eqn:OptimizationProblemObjectiveFunctional} is to introduce additional geometrical constraints, cf. \cite{Schulz2016}. For the investigations in this publication we add a volume and barycenter equality constraint for the shape~$\shape$, i.e.,
\begin{align*}
	\vol(\shape) = \mathcal{V} = \const \quad \text{and} \quad \bary(\shape) = \bm{\mathcal{B}} = \const
\end{align*}
An additional perimeter regularization is commonly added to optimization problems to overcome ill-posedness of the problem, cf.~\cite{LoayzaGeiersbachWelkerHandbook,Schulz2015,Schulz2016a,Sturm2016}. As we already have the possibility to enforce geometrical constraints, we instead include a perimeter constraint
\begin{align*}
	\peri(\shape) = \mathcal{P} = \const
\end{align*}
The augmented Lagrange functional, which includes the geometrical equality constraints enforced by an augmented Lagrange approach with the Lagrange multiplier~$\bm{\lambda}=\left(\lambda_1, \lambda_2, \lambda_3, \lambda_4\right)^\top$ and the penalty parameter $\nu>0$, is defined as
\begin{align}
	\label{eqn:augmentedLagrangeFunctional}
	L_A(\shape,(\velocity,\pressure),(\perturbvelocitytwo,\perturbpressuretwo)) = J(\shape,\velocity) + R(\shape,(\velocity,\pressure),(\perturbvelocitytwo,\perturbpressuretwo))  - \bm{\lambda}^\top \bm{c}(\shape) + \frac{\nu}{2} \left\| \bm{c}(\shape) \right\|_2^2,
\end{align}
where $\bm{c}(\shape) = (\vol(\shape) - \mathcal{V}, \bary(\shape) - \bm{\mathcal{B}} , \peri(\shape) - \mathcal{P} )^\top \in \mathbb{R}^4$ is a vector containing the geometrical constraints.

\section{Optimization technique}
\label{sec:ApplicationSemiderivativeApproach}

It is obvious that \eqref{eqn:regularizedNonDifferentiability} and therefore also \eqref{eqn:augmentedLagrangeFunctional} is non-differentiable in the classical sense due to the presence of the $\max$-operator. However, there is more information available for the optimization algorithm than only using function evaluations, which leads to the concept of Gâteaux semiderivatives and semismoothness (cf. \cite{Hintermueller2002,Hintermueller2006,Stadler2004}).
First, we consider the state equation \eqref{eqn:OptimizationProblemPDEConstraint} of our model and discuss an approach to solve it numerically (cf. \Cref{sec:NumState}).
Then, we focus on the necessary optimality conditions for our non-smooth shape optimization model (cf. \Cref{sec:GateauxOptimalitySystem}), which leads to an optimization technique that can handle the non-differentiability
 (cf. \Cref{sec:OptimizationAlgorithm}). 
 
In the following, we will focus on Gâteaux semiderivatives, e.g., to propose an approach to solve the state equation in \Cref{sec:NumState} and to derive the  %
 adjoint system in \Cref{sec:GAdj}. 
 Moreover, we need the Eulerian derivative, e.g., %
 to set up the design equation in \Cref{sec:DesignEqu}. 
 For the detailed definitions, we refer to the literature, e.g., \cite[Definition~3.1]{Delfour2020} for Gâteaux semiderivative %
 and \cite[Definition~2.19]{SokoZol} for the Eulerian derivative. 
Given a functional $\mathcal{F}(\shape,\bm{x},\bm{z})$ depending on a shape $\shape$ and elements $\bm{x},\bm{z}$ of topological spaces, 
 we will denote the (total) Gâteaux semiderivative of  $\mathcal{F}$ by $d^G\mathcal{F}$.
 The notation $d^G_{\bm{x}}\mathcal{F}(\shape,\bm{x},\bm{z})[\hat{\bm{x}}]$ (and $ d^G_{\bm{z}}\mathcal{F}(\shape,\bm{x},\bm{z})[\hat{\bm{z}}]$) means the Gâteaux semiderivative of  $\mathcal{F}$ with respect to $\bm{x}$ (and $\bm{z}$) in direction $\hat{\bm{x}}$ (and $\hat{\bm{z}}$).
The Eulerian derivative of  $\mathcal{F}$ at $\shape$ in direction $\bm{V}$ is denoted by $d_E \mathcal{F}(\shape,\bm{x},\bm{z})[\bm{V}]$. %
 In the case of Fréchet differentiable terms (for a definition, cf., e.g., \cite[Definition~2.3]{Delfour2020}), the Gâteaux semiderivatives correspond to
  the Fréchet operator applied to the corresponding direction. Moreover, the Eulerian derivative 
  corresponds to the classical shape derivative in case of shape differentiable terms.
 
 \subsection{Newton-like approach for solving the state equation}
 \label{sec:NumState}

In this subsection, we focus on an approach to numerically solve
 \begin{align}
 \label{eqn:OptimizationProblemPDEConstraintOptimalityCondition}
 \begin{aligned}
 0 %
 & = \int_{\compDomain} 2 \viscosity \, \strain(\velocity) : \strain(\adjvelocity) +\density ((\velocity \cdot \nabla) \velocity) \cdot \adjvelocity - \pressure \Div(\adjvelocity)  \\
 & \phantom{=\int_{\compDomain}}+ \adjpressure \Div(\velocity)- \bm{f} \cdot \adjvelocity + \regNonDifferentiability(\strain(\velocity)) : \strain(\adjvelocity) \dx%
 \end{aligned}
 \end{align}
 for the state $(\velocity,\pressure) \in \mathcal{H} \times L^2(\compDomain)$ for all $\adjvelocity, \adjpressure \in \mathcal{H}_0 \times L^2(\compDomain)$. Due to the nonlinear occurrence of the velocity an iterative scheme can be employed for the solution. In the following, we will derive a Gâteaux semiderivative of \eqref{eqn:OptimizationProblemPDEConstraintOptimalityCondition} with respect to $(\velocity,\pressure)$ in the direction $(\perturbvelocity, \perturbpressure) \in \mathcal{H}_0 \times L^2(\compDomain)$, which is required for the formulation of a Newton-like scheme to solve this nonlinear and in the classical sense non-differentiable PDE. If a term is not dependent on $\pressure$, then $d^G_{(\velocity,\pressure)} = d^G_{\velocity}$.
 
 As nearly all terms are differentiable the only difficulty lies in the Gâteaux semiderivative of  $\regNonDifferentiability(\strain(\velocity))$ for which the derivation will be performed in the following. Similar to \cite{Reyes2023} we split the derivation into three cases: %
 \begin{itemize}
\item  \emph{Case 1: $\yieldStress > \regParam \! \left\| \strain(\velocity) \right\|$.} \\We immediately obtain $\max(\yieldStress, \regParam\! \left\| \strain(\velocity) \right\|) = \yieldStress$. This gives the Gâteaux semiderivative
\begin{align*}
d^G_{\velocity} \left( \regParam \, \strain(\velocity) \right) [\perturbvelocity] = \regParam \, \strain(\perturbvelocity).
\end{align*}
\item  \emph{Case 2: $\yieldStress < \regParam \! \left\| \strain(\velocity) \right\|$.} \\
Opposed to the previous case, now $\max(\yieldStress, \regParam\! \left\| \strain(\velocity) \right\|) = \regParam\! \left\| \strain(\velocity) \right\|$, which yields
\begin{align*}
d^G_{\velocity} \left( \yieldStress \frac{\strain(\velocity)}{\left\| \strain(\velocity) \right\|} \right) [\perturbvelocity]
= \frac{\yieldStress \, d^G_{\velocity} (\strain(\velocity))[\perturbvelocity]}{\left\| \strain(\velocity) \right\|} - \frac{\yieldStress \, d^G_{\velocity} (\left\| \strain(\velocity) \right\|)[\perturbvelocity] \, \strain(\velocity)}{\left\| \strain(\velocity) \right\|^2}.
\end{align*}
The Gâteaux semiderivative $d^G_{\velocity} (\left\| \strain(\velocity) \right\|)[\perturbvelocity]$ can be found for example in \cite[Equation~(3.5)]{Delfour2020}, and---adapted to our setting---reads
\begin{align*}
d^G_{\velocity} (\left\| \strain(\velocity) \right\|)[\perturbvelocity] = \begin{cases}
\frac{\strain(\velocity) : \strain(\perturbvelocity)}{\left\| \strain(\velocity) \right\|} &\text{if } \strain(\velocity) \neq \bm{0}, \\
\left\| \strain(\perturbvelocity) \right\| &\text{if } \strain(\velocity) = \bm{0}.
\end{cases}
\end{align*}
However, since $\yieldStress>0$ and therefore $\left\| \strain(\perturbvelocity) \right\| > 0$ in this case, the term for $\strain(\velocity) = \bm{0}$ is not relevant and we get
\begin{align*}
d^G_{\velocity} \left( \yieldStress \frac{\strain(\velocity)}{\left\| \strain(\velocity) \right\|} \right) [\perturbvelocity]
= \yieldStress \frac{\strain(\perturbvelocity) }{\left\| \strain(\velocity) \right\|} - \yieldStress \frac{(\strain(\velocity) : \strain(\perturbvelocity)) \, \strain(\velocity)}{\left\| \strain(\velocity) \right\|^3}.
\end{align*}
\item  \emph{Case 3: $\yieldStress = \regParam \! \left\| \strain(\velocity) \right\|$.}\\
As in \cite[Section~2.1]{GoldammerSchulzWelker} we obtain
\begin{align*}
&d^G_{\velocity} \left( \frac{\yieldStress \regParam \, \strain(\velocity)}{\max(\yieldStress, \regParam\! \left\| \strain(\velocity) \right\|)} \right) [\perturbvelocity] \\
&= \frac{\yieldStress \regParam \, d^G_{\velocity} (\strain(\velocity))[\perturbvelocity]}{\max(\yieldStress, \regParam\! \left\| \strain(\velocity) \right\|)} - \frac{\yieldStress \regParam \, d^G_{\velocity} (\max(\yieldStress, \regParam\! \left\| \strain(\velocity) \right\|))[\perturbvelocity] \, \strain(\velocity)}{\max(\yieldStress, \regParam\! \left\| \strain(\velocity) \right\|)^2} \\
&= \frac{\yieldStress \regParam \, \strain(\perturbvelocity)}{\max(\yieldStress, \regParam\! \left\| \strain(\velocity) \right\|)} - \frac{\yieldStress \regParam \, \max\!\left(0, \regParam \frac{\strain(\velocity) : \strain(\perturbvelocity)}{\left\| \strain(\velocity) \right\|}\right) \, \strain(\velocity)}{\max(\yieldStress, \regParam\! \left\| \strain(\velocity) \right\|)^2} \\
&= \frac{\yieldStress \regParam \, \strain(\perturbvelocity)}{\max(\yieldStress, \regParam\! \left\| \strain(\velocity) \right\|)} - \frac{\yieldStress \regParam^2 \, \max\!\left(0, \strain(\velocity) : \strain(\perturbvelocity)\right) \, \strain(\velocity)}{\max(\yieldStress, \regParam\! \left\| \strain(\velocity) \right\|)^2 \left\| \strain(\velocity) \right\|},
\end{align*}
where the Gâteaux semiderivative of $\max$ was used in the third line, which is %
\begin{align*}
d_x^G \max \! \left(0, x\right)[v] = \begin{cases}
0 & \text{if } x < 0, \\
\max(0,v) & \text{if } x = 0, \\
v & \text{if } x > 0
\end{cases}
\end{align*}
(cf., e.g., \cite{GoldammerSchulzWelker}).
Since $\yieldStress = \regParam \! \left\| \strain(\velocity) \right\|$ we replace the $\max$-terms in the denominators, which in this case results in
\begin{align*}
d^G_{\velocity} \left( \frac{\yieldStress \regParam \, \strain(\velocity)}{\max(\yieldStress, \regParam\! \left\| \strain(\velocity) \right\|)} \right) [\perturbvelocity]
= \yieldStress \frac{\strain(\perturbvelocity)}{\left\| \strain(\velocity) \right\|}  - \yieldStress \frac{\max(0, \strain(\velocity) : \strain(\perturbvelocity)) \, \strain(\velocity)}{\left\| \strain(\velocity) \right\|^3}.
\end{align*}
Therefore, %
\begin{align}
\label{eqn:GateauxDerivativeBingham}
\begin{aligned}
&d^G_{\velocity} \left( \regNonDifferentiability(\strain(\velocity))\right) [\perturbvelocity] 
= \begin{cases}
\regParam \, \strain(\perturbvelocity) : \strain(\adjvelocity) & \text{if } \yieldStress > \regParam \! \left\| \strain(\velocity) \right\|\\
\yieldStress \frac{\strain(\perturbvelocity)}{\left\| \strain(\velocity) \right\|}  - \yieldStress\frac{\max(0, \strain(\velocity) : \strain(\perturbvelocity)) \, \strain(\velocity)}{\left\| \strain(\velocity) \right\|^3} & \text{if } \yieldStress = \regParam \! \left\| \strain(\velocity) \right\|\\
\yieldStress\frac{\strain(\perturbvelocity)}{\left\| \strain(\velocity) \right\|} - \yieldStress\frac{(\strain(\velocity) : \strain(\perturbvelocity)) \, \strain(\velocity)}{\left\| \strain(\velocity) \right\|^3}  & \text{if } \yieldStress < \regParam \! \left\| \strain(\velocity) \right\|
\end{cases} \\
&= \frac{\yieldStress \regParam \, \strain(\perturbvelocity)}{\max(\yieldStress, \regParam\! \left\| \strain(\velocity) \right\|)} - \begin{cases}
\bm{0} & \text{if } \yieldStress > \regParam \! \left\| \strain(\velocity) \right\|, \\
\yieldStress\frac{\max(0, \strain(\velocity) : \strain(\perturbvelocity)) \, \strain(\velocity)}{\left\| \strain(\velocity) \right\|^3} & \text{if } \yieldStress = \regParam \! \left\| \strain(\velocity) \right\|, \\
\yieldStress\frac{(\strain(\velocity) : \strain(\perturbvelocity)) \, \strain(\velocity)}{\left\| \strain(\velocity) \right\|^3} & \text{if } \yieldStress < \regParam \! \left\| \strain(\velocity) \right\|.
\end{cases}
\end{aligned}
\end{align}
 \end{itemize}
 
\smallskip
 By defining the functional
 \begin{align}
 \label{eqn:NonlinearTermStateEquation}
 \bm{d}(\strain(\velocity), \strain(\perturbvelocity)) = \begin{cases}
 \bm{0} & \text{if } \yieldStress > \regParam \! \left\| \strain(\velocity) \right\|, \\
 \yieldStress\frac{\max(0, \strain(\velocity) : \strain(\perturbvelocity)) \, \strain(\velocity)}{\left\| \strain(\velocity) \right\|^3} & \text{if } \yieldStress = \regParam \! \left\| \strain(\velocity) \right\|, \\
 \yieldStress\frac{(\strain(\velocity) : \strain(\perturbvelocity)) \, \strain(\velocity)}{\left\| \strain(\velocity) \right\|^3} & \text{if } \yieldStress < \regParam \! \left\| \strain(\velocity) \right\|
 \end{cases}
 \end{align}
 the Gâteaux semiderivative of \eqref{eqn:OptimizationProblemPDEConstraintOptimalityCondition} with respect to $(\velocity,\pressure)$ in the direction $(\perturbvelocity, \perturbpressure)$ using \eqref{eqn:GateauxDerivativeBingham} and \eqref{eqn:NonlinearTermStateEquation} is given by
 \begin{align}
 \label{eqn:GateauxDerivativeResidual}
 \begin{aligned}
 &d^G_{(\velocity,\pressure)} R(\shape,(\velocity, \pressure),( \adjvelocity, \adjpressure))[(\perturbvelocity, \perturbpressure)] \\
 &= \int_{\compDomain} 2 \viscosity \, \strain(\perturbvelocity) : \strain(\adjvelocity) + ((\perturbvelocity \cdot \nabla) \velocity) \cdot \adjvelocity + \density \,((\velocity \cdot \nabla) \perturbvelocity) \cdot \adjvelocity - \perturbpressure \Div(\adjvelocity) \\
 &\hphantom{= \int_{\compDomain} \ } + \adjpressure \Div(\perturbvelocity) + \frac{\yieldStress \regParam \, \strain(\perturbvelocity) : \strain(\adjvelocity)}{\max(\yieldStress, \regParam\! \left\| \strain(\velocity) \right\|)} - \bm{d}(\strain(\velocity), \strain(\perturbvelocity)) : \strain(\adjvelocity) \dx.
 \end{aligned}
 \end{align}

 A Newton-like method to solve \eqref{eqn:OptimizationProblemPDEConstraintOptimalityCondition}  using the Gâteaux semiderivative is given by
 \begin{align*}
 d^G_{(\velocity,\pressure)} R(\shape, (\velocity^k, \pressure^k), (\adjvelocity, \adjpressure)) &\left[(\perturbvelocity^{k+1},\perturbpressure^{k+1})\right] \\ &= -R(\shape, (\velocity^k, \pressure^k), (\adjvelocity, \adjpressure)) \quad \forall \adjvelocity, \adjpressure \in \mathcal{H}_0 \times L^2(\compDomain) \\
 \left( \velocity^{k+1}, \pressure^{k+1} \right) &= \left( \velocity^k, \pressure^k \right) + \alpha^{k+1} \left( \perturbvelocity^{k+1},\perturbpressure^{k+1} \right)
 \end{align*}
 with an appropriate step size $0<\alpha^{k+1} \leq 1$. For a linear system of equations to compute the update $(\perturbvelocity^{k+1},\perturbpressure^{k+1})$ we require linearity in $(\perturbvelocity,\perturbpressure)$. Thus, we replace the term $\max(0, \strain(\velocity^k) : \strain(\perturbvelocity^{k+1}))$ by $\max(0, \strain(\velocity^k) : \strain(\perturbvelocity^{k}))$, which is reasonable since $(\perturbvelocity^{k},\perturbpressure^{k}) \to \bm{0}$ as $R$ goes to zero with increasing $k$. With the modification  of \eqref{eqn:NonlinearTermStateEquation} we define
 \begin{align*}
 \bm{d}_{mod}(\strain(\velocity^k), \strain(\perturbvelocity^{k}), \strain(\perturbvelocity^{k+1})) = \begin{cases}
 \bm{0} & \text{if } \yieldStress > \regParam \! \left\| \strain(\velocity^k) \right\|, \\
 \yieldStress\frac{\max(0, \strain(\velocity^k) : \strain(\perturbvelocity^k)) \, \strain(\velocity^k)}{\left\| \strain(\velocity^k) \right\|^3} 
 & \text{if } \yieldStress = \regParam \! \left\| \strain(\velocity^k) \right\|, \\
 \yieldStress\frac{(\strain(\velocity^k) : \strain(\perturbvelocity^{k+1})) \, \strain(\velocity^k)}{\left\| \strain(\velocity^k) \right\|^3} & \text{if } \yieldStress < \regParam \! \left\| \strain(\velocity^k) \right\|
 \end{cases}
 \end{align*}
 and propose the method to solve for $(\perturbvelocity^{k+1},\perturbpressure^{k+1})$ as
 \begin{align}
 \label{eqn:SemismoothNewtonStateEquation}
 \begin{aligned}
 &\int_{\compDomain} 2 \viscosity \, \strain(\perturbvelocity^{k+1}) : \strain(\adjvelocity) + \density \, ((\perturbvelocity^{k+1} \cdot \nabla) \velocity^k) \cdot \adjvelocity + \density \, ((\velocity^k \cdot \nabla) \perturbvelocity^{k+1}) \cdot \adjvelocity - \perturbpressure^{k+1} \Div(\adjvelocity) \\
 &\hphantom{\int_{\compDomain} \ } + \adjpressure \Div(\perturbvelocity^{k+1}) + \frac{\yieldStress \regParam \, \strain(\perturbvelocity^{k+1}) : \strain(\adjvelocity)}{\max(\yieldStress, \regParam\! \left\| \strain(\velocity^k) \right\|)} - \bm{d}_{mod}(\strain(\velocity^k), \strain(\perturbvelocity^{k}), \strain(\perturbvelocity^{k+1})) : \strain(\adjvelocity) \dx \\
 &= \int_{\compDomain} 2 \viscosity \, \strain(\velocity^k) : \strain(\adjvelocity) + \density \, ((\velocity^k \cdot \nabla) \velocity^k) \cdot \adjvelocity - \pressure^k \Div(\adjvelocity) + \adjpressure \Div(\velocity^k) %
 - \bm{f} \cdot \adjvelocity \\
 &\hphantom{= \int_{\compDomain} \ }+ \frac{\yieldStress \regParam \, \strain(\velocity^k) : \strain(\adjvelocity)}{\max(\yieldStress, \regParam\! \left\| \strain(\velocity^k) \right\|)} \dx  \quad \forall \adjvelocity, \adjpressure \in \mathcal{H}_0 \times L^2(\compDomain)
 \end{aligned}
 \end{align}
 and then perform the step as
 \begin{align*}
 \left( \velocity^{k+1}, \pressure^{k+1} \right) &= \left( \velocity^k, \pressure^k \right) + \alpha^{k+1} \left( \perturbvelocity^{k+1},\perturbpressure^{k+1} \right),
 \end{align*}
 where the step size $\alpha^{k+1}$ is the largest step size that corresponds to Armijo backtracking 
 using the control parameter $\beta=10^{-4}$, cf. \cite[Chapter 6]{Dennis1996}. The step sizes start at the full Newton step ($\alpha^{k+1} = 1$) and are halved in case the sufficient decrease condition is not fulfilled. We use an incremental stopping criterion for the iteration over $k$. The iterations are stopped when $\left\| \begin{pmatrix}
 \perturbvelocity^{k+1},
 \perturbpressure^{k+1}
 \end{pmatrix} \right\|_1 < \epsilon$ for an $\epsilon \ll 1$. We chose to solve the state equation up to an accuracy of $\epsilon = 10^{-6}$.

\subsection{Optimality system}
\label{sec:GateauxOptimalitySystem}

We derive the optimality system in the approach usually employed for Fr\'echet differentiable problems \cite[Section~2.2]{Borzi-Schulz-2012}. We consider the variables $(\velocity,\pressure)$  as dependent on the shape $u$. Then there holds for feasible states $(\velocity(\shape),\pressure(\shape))$
\begin{align*}
J(\shape,\velocity(\shape)) =J(\shape,(\velocity(\shape),\pressure(\shape)))\coloneqq J(\shape,\velocity(\shape)) + R(\shape,(\velocity(\shape),\pressure(\shape)),(\perturbvelocitytwo,\perturbpressuretwo)) %
\end{align*} 
for all $(\perturbvelocitytwo,\perturbpressuretwo)$. 
Thus, for arbitrary $(\perturbvelocitytwo,\perturbpressuretwo)$ the total Gâteaux semiderivative of the objective $J$
is given by the following expressions, where we suppress explicit arguments $\shape$ in the notation for the sake of compact equations:
\begin{align*}
d^G J(\shape,\velocity)[\perturbshape] &=d^G J(\shape,(\velocity,\pressure))[\perturbshape,(\perturbvelocity,\perturbpressure)] \\
&=d_E J(\shape,\velocity)[\perturbshape] + d^G_{\velocity} J(\shape,\velocity)[d_E \velocity[\perturbshape]] + d_E R(\shape,(\velocity,\pressure),(\perturbvelocitytwo,\perturbpressuretwo))[\perturbshape] \\
&\hphantom{=\,} +d^G_{(\velocity,\pressure)}R(\shape,(\velocity,\pressure),(\perturbvelocitytwo,\perturbpressuretwo))[d_E (\velocity,\pressure)[\perturbshape]],
\end{align*} 
where $\perturbvelocity=d_E \velocity[\perturbshape]$ and $\perturbpressure=d_E \pressure[\perturbshape] $ due to the dependence of $\velocity$ and $\pressure$ on $\shape$. 
If we assume that the shape $\shape$ is optimal, then all directional derivatives are nonnegative, i.e.
\[
d^G J(\shape,\velocity(\shape))[\perturbshape] \ge 0
\] 
for all feasible $\perturbshape$. 

If we apply the formulated system to our model problem, we see 
 that the derivative $d^G_{(\velocity,\pressure)}R(\shape,(\velocity,\pressure),(\perturbvelocitytwo,\perturbpressuretwo))[(\perturbvelocity, \perturbpressure)]$ has a component which is nonlinear in perturbations $(\perturbvelocity, \perturbpressure)$, cf. \eqref{eqn:NonlinearTermStateEquation}--\eqref{eqn:GateauxDerivativeResidual}. This component holds only on the set $\{\velocity\,:\,\yieldStress = \regParam \! \left\| \strain(\velocity)\right\|\}$, which is the boundary of the active set. If this set has Lebesgue measure 0, then this nonlinear term can be replaced by zero and the resulting approximation
 $\widetilde{d^G_{(\velocity,\pressure)}R}(\shape,(\velocity,\pressure),(\perturbvelocitytwo,\perturbpressuretwo))[(\perturbvelocity, \perturbpressure)]$ is linear in $(\perturbvelocity, \perturbpressure)$. Applying the standard approach, we find the approximate adjoint variables $(\perturbvelocitytwo,\perturbpressuretwo)$ as solutions of the linear equation
\[
0= d^G_{\velocity} J(\shape,\velocity)[\perturbvelocity] + \widetilde{d^G_{(\velocity,\pressure)}R}(\shape,(\velocity,\pressure),(\perturbvelocitytwo,\perturbpressuretwo))[(\perturbvelocity,\perturbpressure)]
\]
for all $(\perturbvelocity,\perturbpressure)$ and obtain the approximation to the total Gâteaux semiderivative as
\[
d^G J\approx d_E J+d_E R.
\]
Since this approximation may be prone to error, which may prevent descent properties, we employ a safeguard technique based on a mollified problem formulation similar to~\cite{Luft2020}.

In the following, we set up an optimality system for our VI constrained shape optimization model.
It consists of the the state equation, the Gâteaux adjoint equation and the Gâteaux design equation. The state equation is given from the model itself but could also be achieved by Gâteaux semidifferentiating the augmented Lagrange functional \eqref{eqn:augmentedLagrangeFunctional} with respect to $(\adjvelocity,\adjpressure)$. If we Gâteaux semidifferentiate the augmented Lagrange functional \eqref{eqn:augmentedLagrangeFunctional} with respect to $(\velocity,\pressure)$ we obtain the Gâteaux adjoint equation. The derivation can be found in \Cref{sec:GAdj}. The Gâteaux design equation, which we get be calculating the Eulerian derivative of the augmented Lagrange functional \eqref{eqn:augmentedLagrangeFunctional}, is derived in \Cref{sec:DesignEqu}.
The optimality system is then used to formulate an optimization technique in \Cref{sec:OptimizationAlgorithm}, in which we need to incorporate a safeguard technique based on a smooth regularization of the non-smooth terms of the model. Thus, we will also concentrate on the regularized max-operator in \Cref{sec:RegMax}.

\subsubsection{Gâteaux adjoint equation} 
\label{sec:GAdj}

After obtaining $\velocity \in \mathcal{H}$ and $\pressure \in L^2(\compDomain)$ from solving \eqref{eqn:OptimizationProblemPDEConstraintOptimalityCondition}, these can be used in the adjoint equation, which is derived from the augmented Lagrange functional \eqref{eqn:augmentedLagrangeFunctional} by semidifferentiating in the Gâteaux sense with respect to $(\velocity,\pressure)$ in direction $(\perturbvelocity, \perturbpressure)$. This Gâteaux semiderivative of $R(\shape,(\velocity, \pressure),( \adjvelocity, \adjpressure))$ has already been obtained in \eqref{eqn:GateauxDerivativeResidual} and the additional terms from the augmented Lagrange approach do not depend neither on $\velocity$ nor $\pressure$. Thus, we only need to further obtain the Gâteaux semiderivative of $J(\shape,\velocity)$ with respect to $\velocity$ in direction $\perturbvelocity$, which reads
\begin{align*}
	d^G_{\velocity} \left( \int_{\compDomain} \frac{\viscosity}{2} \nabla \velocity : \nabla \velocity \dx \right)[\perturbvelocity] = \int_{\compDomain} \viscosity \, \nabla \perturbvelocity : \nabla \velocity \dx.
\end{align*}
Combining this with the Gâteaux semiderivative of $R$ in \eqref{eqn:GateauxDerivativeResidual} yields the Gâteaux semiderivative of $L_A$ with respect to $(\velocity,\pressure)$ in direction $(\perturbvelocity, \perturbpressure)$:
\begin{align}
	\label{eqn:GateauxAdjointEquation}
	\begin{aligned}
		&d^G_{(\velocity,\pressure)} \left( L_A(\shape, (\velocity,\pressure),(\adjvelocity,\adjpressure)) \right)[(\perturbvelocity, \perturbpressure)] \\
		&= \int_{\compDomain} \viscosity \, \nabla \perturbvelocity : \nabla \velocity \dx + \int_{\compDomain} 2 \viscosity \, \strain(\perturbvelocity) : \strain(\adjvelocity) + \density ((\perturbvelocity \cdot \nabla) \velocity) \cdot \adjvelocity + \density ((\velocity \cdot \nabla) \perturbvelocity) \cdot \adjvelocity  \\
		&\hphantom{=\int_{\compDomain} \viscosity \, \nabla \perturbvelocity : \nabla \velocity \dx + \int_{\compDomain}\ } - \perturbpressure \Div(\adjvelocity) + \adjpressure \Div(\perturbvelocity) + \frac{\yieldStress \regParam \, \strain(\perturbvelocity) : \strain(\adjvelocity)}{\max(\yieldStress, \regParam\! \left\| \strain(\velocity) \right\|)} \\
		&\hphantom{=\int_{\compDomain} \viscosity \, \nabla \perturbvelocity : \nabla \velocity \dx + \int_{\compDomain}\ } - \bm{d}(\strain(\velocity), \strain(\perturbvelocity)) : \strain(\adjvelocity) \dx
	\end{aligned}
\end{align}
This expression is not linear with respect to the direction $(\perturbvelocity, \perturbpressure)$ because of the component $\bm{d}(\strain(\velocity), \strain(\perturbvelocity))$. We replace this by a linear approximation 
\begin{align*}
	\widetilde{\bm{d}}(\strain(\velocity), \strain(\perturbvelocity)) = \begin{cases}
		\bm{0} & \text{if } \yieldStress > \regParam \! \left\| \strain(\velocity) \right\| \\
		\bm{0} & \text{if } \yieldStress = \regParam \! \left\| \strain(\velocity) \right\| \\
		\yieldStress\frac{(\strain(\velocity) : \strain(\perturbvelocity)) \, \strain(\velocity)}{\left\| \strain(\velocity) \right\|^3} & \text{if } \yieldStress < \regParam \! \left\| \strain(\velocity) \right\|
	\end{cases}
\end{align*}
as mentioned at the beginning of \Cref{sec:GateauxOptimalitySystem} in order to solve the linear equation
\begin{align}
	\label{eqn:GateauxAdjointEquationApprox}
	\begin{aligned}
		& \widetilde{d^G_{(\velocity,\pressure)} L_A}(\shape, (\velocity,\pressure),(\adjvelocity,\adjpressure))[(\perturbvelocity, \perturbpressure)] = 0 \\
		&= \int_{\compDomain} \viscosity \, \nabla \perturbvelocity : \nabla \velocity \dx + \int_{\compDomain} 2 \viscosity \, \strain(\perturbvelocity) : \strain(\adjvelocity) + \density ((\perturbvelocity \cdot \nabla) \velocity) \cdot \adjvelocity + \density ((\velocity \cdot \nabla) \perturbvelocity) \cdot \adjvelocity  \\
		&\hphantom{=\int_{\compDomain} \viscosity \, \nabla \perturbvelocity : \nabla \velocity \dx + \int_{\compDomain}\ } - \perturbpressure \Div(\adjvelocity) + \adjpressure \Div(\perturbvelocity) + \frac{\yieldStress \regParam \, \strain(\perturbvelocity) : \strain(\adjvelocity)}{\max(\yieldStress, \regParam\! \left\| \strain(\velocity) \right\|)} \\
		&\hphantom{=\int_{\compDomain} \viscosity \, \nabla \perturbvelocity : \nabla \velocity \dx + \int_{\compDomain}\ } - \widetilde{\bm{d}}(\strain(\velocity), \strain(\perturbvelocity)) : \strain(\adjvelocity) \dx
	\end{aligned}
\end{align}
for all $ \perturbvelocity, \perturbpressure \in \mathcal{H}_0 \times L^2(\compDomain)$ to obtain an approximate adjoint $\adjvelocity,\adjpressure \in \mathcal{H}_0 \times L^2(\compDomain)$.

\subsubsection{Gâteaux design equation} 
\label{sec:DesignEqu}

For optimization with respect to the shape $\shape$ %
the Eulerian derivative %
is required. We use a material derivative approach \cite{Berggren2009} and its extensions in~\cite[Section~2.2.1]{GeiersbachSuchanWelker_NS} for this derivation with the material derivative in direction~$\meshDeformation$ denoted as $d_m(\cdot)[\meshDeformation]$ or $\dot{(\cdot)}$. The interested reader is referred to, e.g., \cite{GeiersbachSuchanWelker_NS} for the derivation of the differentiable terms in the classical sense. The shape derivative of the terms from the augmented Lagrange method has already been derived in, e.g., \cite{GeiersbachSuchanWelker_NS,Schulz2016,SokoZol}. The computation of the Eulerian derivative of the nondifferentiable term $\bm{d}(\strain(\velocity), \strain(\perturbvelocity))$ is provided in the following.

We start with perturbing the functional $\int_\compDomain \regNonDifferentiability(\strain(\velocity)) : \strain(\adjvelocity) \dx$ by $\meshDeformation$ and obtain
\begin{align*}
	&d_E \left( \int_\compDomain \regNonDifferentiability(\strain(\velocity)) : \strain(\adjvelocity) \dx \right) \left[\meshDeformation\right] \\
	&= \int_{\compDomain} d_m \! \left( \regNonDifferentiability(\strain(\velocity)) \right) [\meshDeformation] : \strain(\adjvelocity) + \regNonDifferentiability(\strain(\velocity)) : d_m \! \left( \strain(\adjvelocity) \right) [\meshDeformation] \\
	&\hphantom{= \int_{\compDomain}\ } + \Div(\meshDeformation) \left( \regNonDifferentiability(\strain(\velocity)) : \strain(\adjvelocity) \right) \dx.
\end{align*}
The material derivative $d_m \! \left(\strain(\adjvelocity)\right) [\meshDeformation]$ with the definition of $\strain(\adjvelocity)$ as in~\eqref{eqn:symmetricGradient} reads
\begin{equation}
\begin{split}
	d_m \! \left(\strain(\adjvelocity)\right) [\meshDeformation] & %
	= \frac{1}{2} \left( \nabla \dot{\adjvelocity} - \nabla \adjvelocity \nabla \meshDeformation + {\nabla \dot{\adjvelocity}}^\top - \left( \nabla \adjvelocity \nabla \meshDeformation \right)^\top \right) \\
	&= \strain(\dot{\adjvelocity}) - \frac{1}{2} \left( \nabla \adjvelocity \nabla \meshDeformation + {\nabla \meshDeformation}^\top {\nabla \adjvelocity}^\top \right).
	\label{eqn:MaterialDerivativeSymmetricGradient}
\end{split}
\end{equation}
The material derivative of $\regNonDifferentiability(\strain(\velocity)) = \frac{\yieldStress \regParam \, \strain(\velocity)}{\max(\yieldStress, \regParam\! \left\| \strain(\velocity) \right\|)}$, using the quotient material derivative rule, which follows with a similar derivation as the shape derivative quotient rule in \cite[Section~2.2.1]{GeiersbachSuchanWelker_NS}, reads
\begin{align}
	\label{eqn:MaterialDerivativeHGamma}
	\begin{aligned}
	&d_m \! \left( \frac{\yieldStress \regParam \, \strain(\velocity)}{\max(\yieldStress, \regParam\! \left\| \strain(\velocity) \right\|)} \right) [\meshDeformation] \\
	&= \frac{\yieldStress \regParam \, d_m ( \strain(\velocity) ) [\meshDeformation]}{\max(\yieldStress, \regParam\! \left\| \strain(\velocity) \right\|)} - \frac{\yieldStress \regParam \, \strain(\velocity) \, d_m (\max(\yieldStress, \regParam\! \left\| \strain(\velocity) \right\|)) [\meshDeformation]}{\max(\yieldStress, \regParam\! \left\| \strain(\velocity) \right\|)^2}.
	\end{aligned}
\end{align}
The first term in \eqref{eqn:MaterialDerivativeHGamma} only requires
$$d_m \! \left(\strain(\velocity)\right) [\meshDeformation] = \strain(\dot{\velocity}) - \frac{1}{2} \left( \nabla \velocity \nabla \meshDeformation + {\nabla \meshDeformation}^\top {\nabla \velocity}^\top \right),$$
which follows from \eqref{eqn:MaterialDerivativeSymmetricGradient} directly. For the second term, we use the fact that $\max(\yieldStress, \regParam\! \left\| \strain(\velocity) \right\|) = g + \max(0, \regParam\! \left\| \strain(\velocity) \right\| - \yieldStress)$ and obtain
\begin{align}
	\label{eqn:MaterialDerivativeMax}
	\begin{aligned}
		d_m (\max(\yieldStress, \regParam\! \left\| \strain(\velocity) \right\|)) [\meshDeformation] %
		= \begin{cases}
			0 & \text{if } \yieldStress > \regParam \! \left\| \strain(\velocity) \right\|, \\
			\regParam \max(0, d_m ( \left\| \strain(\velocity) \right\| ) [\meshDeformation] ) & \text{if } \yieldStress = \regParam \! \left\| \strain(\velocity) \right\|, \\
			\regParam \, d_m ( \left\| \strain(\velocity) \right\| ) [\meshDeformation] & \text{if } \yieldStress < \regParam \! \left\| \strain(\velocity) \right\|.
		\end{cases}
	\end{aligned}
\end{align}
With the definition of the Frobenius product $\left\| \strain(\velocity) \right\| = \sqrt{ \strain(\velocity) : \strain(\velocity)}$, we get
\begin{equation}
\begin{split}
&	d_m \! \left( \sqrt{ \strain(\velocity) : \strain(\velocity)}\right) [\meshDeformation] 
\\&= \begin{cases}
		2 \strain(\velocity) : \frac{d_m ( \strain(\velocity) ) [\meshDeformation]}{2 \sqrt{ \strain(\velocity) : \strain(\velocity) }} = \frac{\strain(\velocity) : d_m ( \strain(\velocity))[\meshDeformation] }{\left\| \strain(\velocity) \right\|} & \text{if } \left\| \strain(\velocity) \right\| \neq 0, \\
		\| d_m ( \strain(\velocity))[\meshDeformation] \| & \text{if } \left\| \strain(\velocity) \right\| = 0.
	\end{cases}
	\end{split}
\end{equation}
As in \Cref{sec:NumState}, case~3, since $\left\| \strain(\velocity) \right\| \geq g > 0$, this gives
\begin{align}
	\label{eqn:MaterialDerivativeFrobeniusNorm}
	d_m \! \left( \sqrt{ \strain(\velocity) : \strain(\velocity)}\right) [\meshDeformation] = \frac{\strain(\velocity) : d_m ( \strain(\velocity))[\meshDeformation] }{\left\| \strain(\velocity) \right\|}.
\end{align}
Combining \eqref{eqn:MaterialDerivativeHGamma}, \eqref{eqn:MaterialDerivativeMax} and \eqref{eqn:MaterialDerivativeFrobeniusNorm} then gives the material derivative of $\regNonDifferentiability(\strain(\velocity))$ as
\begin{align*}
	\begin{aligned}
	&d_m \! \left( \regNonDifferentiability(\strain(\velocity)) \right) [\meshDeformation] \\
	&= \frac{\yieldStress \regParam \left( \strain(\dot{\velocity}) - \frac{1}{2} \left( \nabla \velocity \nabla \meshDeformation + {\nabla \meshDeformation}^\top {\nabla \velocity}^\top \right) \right)}{\max(\yieldStress, \regParam\! \left\| \strain(\velocity) \right\|)} \\
	&\hphantom{=\ } - \begin{cases}
		\bm{0} & \text{if } \yieldStress > \regParam \! \left\| \strain(\velocity) \right\|, \\
		\frac{\yieldStress \regParam^2 \, \strain(\velocity) \max \left(0, \strain(\velocity) : \left( \strain(\dot{\velocity}) - \frac{1}{2} \left( \nabla \velocity \nabla \meshDeformation + {\nabla \meshDeformation}^\top {\nabla \velocity}^\top \right) \right) \right)}{\max(\yieldStress, \regParam\! \left\| \strain(\velocity) \right\|)^2 \, \left\| \strain(\velocity) \right\|} & \text{if } \yieldStress = \regParam \! \left\| \strain(\velocity) \right\|,  \\
		\frac{\yieldStress \regParam^2 \, \strain(\velocity) \left( \strain(\velocity) : \left( \strain(\dot{\velocity}) - \frac{1}{2} \left( \nabla \velocity \nabla \meshDeformation + {\nabla \meshDeformation}^\top {\nabla \velocity}^\top \right) \right) \right)}{\max(\yieldStress, \regParam\! \left\| \strain(\velocity) \right\|)^2 \, \left\| \strain(\velocity) \right\|} & \text{if } \yieldStress < \regParam \! \left\| \strain(\velocity) \right\|.
	\end{cases}
	\end{aligned}
\end{align*}
Using the fact that $\max(\yieldStress, \regParam\! \left\| \strain(\velocity) \right\| ) = \regParam\! \left\| \strain(\velocity) \right\|$ for $\regParam\! \left\| \strain(\velocity) \right\| \geq g$, we obtain the final result
\begin{align*}
		\begin{aligned}
		&d_m \! \left( \regNonDifferentiability(\strain(\velocity)) \right) [\meshDeformation] \\
		&= \frac{\yieldStress \regParam \left( \strain(\dot{\velocity}) - \frac{1}{2} \left( \nabla \velocity \nabla \meshDeformation + {\nabla \meshDeformation}^\top {\nabla \velocity}^\top \right) \right)}{\max(\yieldStress, \regParam\! \left\| \strain(\velocity) \right\|)} \\
		&\hphantom{=\ } - \begin{cases}
			\bm{0} & \text{if } \yieldStress > \regParam \! \left\| \strain(\velocity) \right\|, \\
			\frac{\yieldStress \, \strain(\velocity) \max \left(0, \strain(\velocity) : \left( \strain(\dot{\velocity}) - \frac{1}{2} \left( \nabla \velocity \nabla \meshDeformation + {\nabla \meshDeformation}^\top {\nabla \velocity}^\top \right) \right) \right)}{\left\| \strain(\velocity) \right\|^3} & \text{if } \yieldStress = \regParam \! \left\| \strain(\velocity) \right\|,  \\
			\frac{\yieldStress \, \strain(\velocity) \left( \strain(\velocity) : \left( \strain(\dot{\velocity}) - \frac{1}{2} \left( \nabla \velocity \nabla \meshDeformation + {\nabla \meshDeformation}^\top {\nabla \velocity}^\top \right) \right) \right)}{\left\| \strain(\velocity) \right\|^3} & \text{if } \yieldStress < \regParam \! \left\| \strain(\velocity) \right\|.
		\end{cases}
	\end{aligned}
\end{align*}
Defining
\begin{align}
	\label{eqn:NonlinearTermDesignEquationM}
	\bm{m} (\velocity, \strain(\dot{\velocity}), \meshDeformation) = \begin{cases}
		\bm{0} & \text{if } \yieldStress > \regParam \! \left\| \strain(\velocity) \right\|, \\
		\frac{\yieldStress \, \strain(\velocity) \max \left(0, \strain(\velocity) : \left( \strain(\dot{\velocity}) - \frac{1}{2} \left( \nabla \velocity \nabla \meshDeformation + {\nabla \meshDeformation}^\top {\nabla \velocity}^\top \right) \right) \right)}{\left\| \strain(\velocity) \right\|^3} & \text{if } \yieldStress = \regParam \! \left\| \strain(\velocity) \right\|,  \\
		\frac{\yieldStress \, \strain(\velocity) \left( \strain(\velocity) : \left( \strain(\dot{\velocity}) - \frac{1}{2} \left( \nabla \velocity \nabla \meshDeformation + {\nabla \meshDeformation}^\top {\nabla \velocity}^\top \right) \right) \right)}{\left\| \strain(\velocity) \right\|^3} & \text{if } \yieldStress < \regParam \! \left\| \strain(\velocity) \right\|
	\end{cases}
\end{align}
and combining the result with all Eulerian derivatives of the differentiable terms then yields the Eulerian derivative of $J(\shape,\velocity)$ in direction $\meshDeformation$
\begin{align*}
	&d_E J(\shape, (\velocity, \pressure)) [\meshDeformation] \\
	&= \int_{\compDomain} 2 \viscosity \, \strain(\velocity) : \strain(\dot{\perturbvelocitytwo}) + \density ((\velocity \cdot \nabla) \velocity) \cdot \dot{\perturbvelocitytwo} - \pressure \Div(\dot{\perturbvelocitytwo}) + \dot{\perturbpressuretwo} \Div(\velocity) - \bm{f} \cdot \dot{\perturbvelocitytwo} \\
	&\hphantom{= \int_{\compDomain}\ }+ \frac{\yieldStress \regParam \, \strain(\velocity) : \strain(\dot{\adjvelocity})}{\max(\yieldStress, \regParam\! \left\| \strain(\velocity) \right\|)} \dx \\
	&\hphantom{=\ } + \int_{\compDomain} \viscosity \, \nabla \dot{\velocity} : \nabla \velocity + 2 \viscosity \, \strain(\dot{\velocity}) : \strain(\adjvelocity) + \density ((\dot{\velocity} \cdot \nabla) \velocity) \cdot \adjvelocity + \density ((\velocity \cdot \nabla) \dot{\velocity}) \cdot \adjvelocity \\
	&\hphantom{=+\int_{\compDomain} \ } - \dot{\pressure} \Div(\adjvelocity) + \adjpressure \Div(\dot{\velocity}) + \frac{\yieldStress \regParam \, \strain(\dot{\velocity}) : \strain(\adjvelocity)}{\max(\yieldStress, \regParam\! \left\| \strain(\velocity) \right\|)} - \bm{m} (\velocity, \strain(\dot{\velocity}), \meshDeformation) : \strain(\adjvelocity) \dx \\
	&\hphantom{=\ } + \int_{\compDomain} - \viscosity \left( \nabla \velocity \nabla \meshDeformation \right) : \nabla \velocity
	- \viscosity \left( \nabla \velocity \nabla \meshDeformation + {\nabla \meshDeformation}^\top {\nabla \velocity}^\top  \right) : \strain(\adjvelocity) \\
	&\hphantom{=+\int_{\compDomain} \ } - \viscosity \, \strain(\velocity) : \left( \nabla \adjvelocity \nabla \meshDeformation + {\nabla \meshDeformation}^\top {\nabla \adjvelocity}^\top  \right) \\
	&\hphantom{=+\int_{\compDomain} \,} - \frac{\yieldStress \regParam \left( \nabla \velocity \nabla \meshDeformation + {\nabla \meshDeformation}^\top {\nabla \velocity}^\top  \right) : \strain(\adjvelocity)}{2\max(\yieldStress, \regParam\! \left\| \strain(\velocity) \right\|)} - \frac{\yieldStress \regParam \, \strain(\velocity) : \left( \nabla \adjvelocity \nabla \meshDeformation + {\nabla \meshDeformation}^\top {\nabla \adjvelocity}^\top  \right)}{2\max(\yieldStress, \regParam\! \left\| \strain(\velocity) \right\|)} \\
	&\hphantom{=+\int_{\compDomain} \,} - \density \left( \nabla \meshDeformation \velocity \cdot \nabla \right) \velocity \cdot \adjvelocity + \pressure {\nabla \adjvelocity}^\top : \nabla \meshDeformation - \left(\nabla \bm{f} \, \meshDeformation \right) \cdot \adjvelocity - \adjpressure {\nabla \velocity}^\top : \nabla \meshDeformation \\
	&\hphantom{=+\int_{\compDomain} \,} + \Div(\meshDeformation) \left( \frac{\viscosity}{2} \nabla \velocity : \nabla \velocity + 2 \viscosity \, \strain(\velocity) : \strain(\adjvelocity)  + \frac{\yieldStress \regParam \, \strain(\velocity) : \strain(\adjvelocity)}{\max(\yieldStress, \regParam\! \left\| \strain(\velocity) \right\|)} \right. \\
	&\hphantom{=+\int_{\compDomain} + \Div(\meshDeformation) \ \ }  \left. \vphantom{\frac{\strain(\adjvelocity)}{\max(\left\| \strain(\velocity) \right\|)}} + \density \left(\velocity \cdot \nabla \right) \velocity \cdot \adjvelocity - \pressure \Div(\adjvelocity) - \bm{f} \cdot \adjvelocity + \adjpressure \Div(\velocity) \right) \dx.
\end{align*}
After inserting the weak form of the state equation~\eqref{eqn:OptimizationProblemPDEConstraint} for the first integral and the weak form of the Gâteaux adjoint equation~\eqref{eqn:GateauxAdjointEquation} for the second integral, we obtain
\begin{align*}
	\begin{aligned}
	&d_E J(\shape,(\velocity,\pressure))[\meshDeformation] \\
	&= \int_{\compDomain} - \viscosity \left( \nabla \velocity \nabla \meshDeformation \right) : \nabla \velocity
	- \viscosity \left( \nabla \velocity \nabla \meshDeformation + {\nabla \meshDeformation}^\top {\nabla \velocity}^\top  \right) : \strain(\adjvelocity) \\
	&\hphantom{=\int_{\compDomain} \,} - \viscosity \, \strain(\velocity) : \left( \nabla \adjvelocity \nabla \meshDeformation + {\nabla \meshDeformation}^\top {\nabla \adjvelocity}^\top  \right) \\
	&\hphantom{=\int_{\compDomain} \,} - \frac{\yieldStress \regParam \left( \nabla \velocity \nabla \meshDeformation + {\nabla \meshDeformation}^\top {\nabla \velocity}^\top  \right) : \strain(\adjvelocity)}{2\max(\yieldStress, \regParam\! \left\| \strain(\velocity) \right\|)} - \frac{\yieldStress \regParam \, \strain(\velocity) : \left( \nabla \adjvelocity \nabla \meshDeformation + {\nabla \meshDeformation}^\top {\nabla \adjvelocity}^\top  \right)}{2\max(\yieldStress, \regParam\! \left\| \strain(\velocity) \right\|)} \\
	&\hphantom{=\int_{\compDomain} \,} - \density \left( \nabla \meshDeformation \velocity \cdot \nabla \right) \velocity \cdot \adjvelocity + \pressure {\nabla \adjvelocity}^\top : \nabla \meshDeformation - \left(\nabla \bm{f} \, \meshDeformation \right) \cdot \adjvelocity - \adjpressure {\nabla \velocity}^\top : \nabla \meshDeformation \\
	&\hphantom{=\int_{\compDomain} \,} + \Div(\meshDeformation) \left( \frac{\viscosity}{2} \nabla \velocity : \nabla \velocity + 2 \viscosity \, \strain(\velocity) : \strain(\adjvelocity)  + \frac{\yieldStress \regParam \, \strain(\velocity) : \strain(\adjvelocity)}{\max(\yieldStress, \regParam\! \left\| \strain(\velocity) \right\|)} \right. \\
	&\hphantom{=\int_{\compDomain} + \Div(\meshDeformation) \ \ \ }  \left. \vphantom{\frac{\strain(\adjvelocity)}{\max(\left\| \strain(\velocity) \right\|)}} + \density \left(\velocity \cdot \nabla \right) \velocity \cdot \adjvelocity - \pressure \Div(\adjvelocity) - \bm{f} \cdot \adjvelocity + \adjpressure \Div(\velocity) \right) \\
	&\hphantom{= \int_{\compDomain}\ } + \left( \bm{d}(\strain(\velocity), \strain(\dot{\velocity})) - \bm{m} (\velocity, \strain(\dot{\velocity}), \meshDeformation) \right) : \strain(\adjvelocity) \dx.%
	\end{aligned}
\end{align*}
The last term is nonlinear in $\meshDeformation$ and still contains material derivative terms that do not cancel out:
\begin{align}
	\label{eqn:NonlinearTermDesignEquationQ}
	\begin{aligned}
	&\bm{q} (\velocity, \strain(\dot{\velocity}), \meshDeformation) \coloneqq \bm{d}(\strain(\velocity), \strain(\dot{\velocity})) - \bm{m} (\velocity, \strain(\dot{\velocity}), \meshDeformation) \\
	&= \begin{cases}
		\bm{0} & \text{if } \yieldStress > \regParam \! \left\| \strain(\velocity) \right\|, \\
		\begin{aligned}[c]
			&\yieldStress\frac{\strain(\velocity)}{\left\| \strain(\velocity) \right\|^3} \Big(\max \left(0, \strain(\velocity) : \strain(\dot{\velocity})\right) \\
			&- \max\left( 0, \strain(\velocity) : \left(\strain(\dot{\velocity}) - \tfrac{1}{2}\! \left( \nabla \velocity \nabla \meshDeformation + {\nabla \meshDeformation}^\top {\nabla \velocity}^\top  \right) \right) \right) \Big)
		\end{aligned} & \text{if } \yieldStress = \regParam \! \left\| \strain(\velocity) \right\|, \\
	 \yieldStress\frac{\strain(\velocity)}{2 \left\| \strain(\velocity) \right\|^3} \, \strain(\velocity) : \left( \nabla \velocity \nabla \meshDeformation + {\nabla \meshDeformation}^\top {\nabla \velocity}^\top  \right) & \text{if } \yieldStress < \regParam \! \left\| \strain(\velocity) \right\|.
	\end{cases}
	\end{aligned}
\end{align}
As in \Cref{sec:GAdj} we replace the nonlinear term $\bm{q} (\strain(\velocity), \strain(\dot{\velocity}), \meshDeformation)$ by a linear approximation
\begin{align*}
	\widetilde{\bm{q}} (\velocity, \meshDeformation)
	= \begin{cases}
		\bm{0} & \text{if } \yieldStress > \regParam \! \left\| \strain(\velocity) \right\|, \\
		\bm{0} & \text{if } \yieldStress = \regParam \! \left\| \strain(\velocity) \right\|, \\
		\yieldStress\frac{\strain(\velocity)}{2 \left\| \strain(\velocity) \right\|^3} \, \strain(\velocity) : \left( \nabla \velocity \nabla \meshDeformation + {\nabla \meshDeformation}^\top {\nabla \velocity}^\top  \right) & \text{if } \yieldStress < \regParam \! \left\| \strain(\velocity) \right\|
	\end{cases}
\end{align*}
to obtain the approximated Eulerian derivative of the augmented Lagrange functional in the direction of a sufficiently smooth vector field $\meshDeformation$ with the unit outward normal $\bm{n}$ as
\begin{align}
	\label{eqn:GateauxShapeDerivativeApprox}
	\begin{aligned}
			&\widetilde{d_E L_A}(\shape,(\velocity,\pressure),(\adjvelocity,\adjpressure))[\meshDeformation] \\
			&= \int_{\compDomain} - \viscosity \left( \nabla \velocity \nabla \meshDeformation \right) : \nabla \velocity
			- \viscosity \left( \nabla \velocity \nabla \meshDeformation + {\nabla \meshDeformation}^\top {\nabla \velocity}^\top  \right) : \strain(\adjvelocity) \\
			&\hphantom{=\int_{\compDomain} \,} - \viscosity \, \strain(\velocity) : \left( \nabla \adjvelocity \nabla \meshDeformation + {\nabla \meshDeformation}^\top {\nabla \adjvelocity}^\top  \right) + \widetilde{\bm{q}} (\velocity, \meshDeformation) : \strain(\adjvelocity) \\
			&\hphantom{=\int_{\compDomain} \,} - \frac{\yieldStress \regParam \left( \nabla \velocity \nabla \meshDeformation + {\nabla \meshDeformation}^\top {\nabla \velocity}^\top  \right) : \strain(\adjvelocity)}{2\max(\yieldStress, \regParam\! \left\| \strain(\velocity) \right\|)} - \frac{\yieldStress \regParam \, \strain(\velocity) : \left( \nabla \adjvelocity \nabla \meshDeformation + {\nabla \meshDeformation}^\top {\nabla \adjvelocity}^\top  \right)}{2\max(\yieldStress, \regParam\! \left\| \strain(\velocity) \right\|)} \\
			&\hphantom{=\int_{\compDomain} \,} - \density \left( \nabla \meshDeformation \velocity \cdot \nabla \right) \velocity \cdot \adjvelocity + \pressure {\nabla \adjvelocity}^\top : \nabla \meshDeformation - \left(\nabla \bm{f} \, \meshDeformation \right) \cdot \adjvelocity - \adjpressure {\nabla \velocity}^\top : \nabla \meshDeformation \\
			&\hphantom{=\int_{\compDomain} \,} + \Div(\meshDeformation) \left( \frac{\viscosity}{2} \nabla \velocity : \nabla \velocity + 2 \viscosity \, \strain(\velocity) : \strain(\adjvelocity)  + \frac{\yieldStress \regParam \, \strain(\velocity) : \strain(\adjvelocity)}{\max(\yieldStress, \regParam\! \left\| \strain(\velocity) \right\|)} \right. \\
			&\hphantom{=\int_{\compDomain} + \Div(\meshDeformation) \ \ }  \left. \vphantom{\frac{\strain(\adjvelocity)}{\max(\left\| \strain(\velocity) \right\|)}} + \density \left(\velocity \cdot \nabla \right) \velocity \cdot \adjvelocity - \pressure \Div(\adjvelocity) - \bm{f} \cdot \adjvelocity + \adjpressure \Div(\velocity) \right) \dx \\
			&\hphantom{= \ } + \left( \lambda_1 - \nu \left(\vol(\shape) - \mathcal{V} \right) \right) \int_{\shape} \meshDeformation^\top \bm{n} \,\mathrm{d}s \\
			&\hphantom{= \ } + \frac{1}{\vol(\shape)} \left( \nu \left(\bary(\shape) - \bm{\mathcal{B}})^\top  - \left(\lambda_2, \lambda_3 \right) \right) \right) \int_{\shape} (\bary(\shape)-\bm{x}) \, \meshDeformation^\top \bm{n} \,\mathrm{d}s \\
			&\hphantom{= \ } + \left( \lambda_4 - \nu \left(\peri(\shape) - \mathcal{P} \right) \right) \int_{\shape} \bm{n}^\top {\nabla \meshDeformation} \bm{n} - \Div(\meshDeformation) \,\mathrm{d}s.
		\end{aligned}
\end{align}

\subsubsection{Regularized $\max$-operator}
\label{sec:RegMax}

As already mentioned in the introduction to \Cref{sec:GateauxOptimalitySystem} we would like to additionally implement a safeguard in case the numerical implementations of the state, adjoint and design equation fail to provide a descent direction for the algorithm on $\admissibleShapeSpace$. Additionally, a comparison between computations using the unregularized approach and using a regularization is another objective of the numerical section. Thus, a regularization of the $\max$-operator is required to obtain a functional that can be Fréchet differentiated.
We use the regularized $\max$-operator in \cite{Reyes2011}, which reads
\begin{align}
	\label{eqn:RegularizationMaxOperator}
	\max\nolimits_{\regParamMax} \left(\yieldStress, \regParam\! \left\| \strain(\velocity) \right\| \right) = \begin{cases}
		\regParam \strain(\velocity) & \text{if } \yieldStress - \frac{1}{2 \regParamMax} \geq \regParam \! \left\| \strain(\velocity) \right\|, \\
		g+\frac{\regParamMax}{2} \left( \regParam \left\| \strain(\velocity) \right\| - g + \frac{1}{2 \regParamMax} \right)^2 & \text{if } \left| \regParam \! \left\| \strain(\velocity) \right\| - \yieldStress \right| \leq \frac{1}{2 \regParamMax}, \\
		\yieldStress \frac{\strain(\velocity)}{\left\| \strain(\velocity) \right\|} & \text{if } \yieldStress + \frac{1}{2 \regParamMax} \leq \regParam \! \left\| \strain(\velocity) \right\|. \\
	\end{cases}
\end{align}

\emph{Regularized state equation.}
Using the regularization in \eqref{eqn:RegularizationMaxOperator} we can determine the 
regularized state equation. By defining
\begin{align*}
	\regNonDifferentiability^\regParamMax(\strain(\velocity)) = \frac{\yieldStress \regParam \, \strain(\velocity)}{\max\nolimits_\regParamMax(\yieldStress, \regParam\! \left\| \strain(\velocity) \right\|)}
\end{align*}
the state equation in the regularized case reads
\begin{align}
	\label{eqn:OptimizationProblemPDEConstraintRegularized}
	\begin{aligned}
		0 = R^\regParamMax(\shape,(\velocity,\pressure),(\adjvelocity,\adjpressure)) = \int_{\compDomain} &2 \viscosity \, \strain(\velocity) : \strain(\adjvelocity) + ((\velocity \cdot \nabla) \velocity) \cdot \adjvelocity - \pressure \Div(\adjvelocity)  \\
		&+ \adjpressure \Div(\velocity) - \bm{f} \cdot \adjvelocity + \regNonDifferentiability^\regParamMax(\strain(\velocity)) : \strain(\adjvelocity) \dx.%
	\end{aligned}
\end{align}
The regularized augmented Lagrange functional then follows as
\begin{align}
	\label{eqn:augmentedLagrangeFunctionalRegularized}
	L_A^\regParamMax(\shape,(\velocity,\pressure),(\perturbvelocitytwo,\perturbpressuretwo)) = J(\shape,\velocity) + R^\regParamMax(\shape,(\velocity,\pressure),(\perturbvelocitytwo,\perturbpressuretwo))  - \bm{\lambda}^\top \bm{c}(\shape) + \frac{\nu}{2} \left\| \bm{c}(\shape) \right\|_2^2.
\end{align}
As the regularization has not removed the nonlinear dependence of $R^\regParamMax$ on $(\velocity,\pressure)$ we require a classical Newton's method, for which the Fréchet operator with respect to $(\velocity,\pressure)$ of $R^\regParamMax(\shape,(\velocity,\pressure),(\adjvelocity,\adjpressure))$ applied to the direction $(\perturbvelocity, \perturbpressure)$ is needed. The derivative of the differentiable terms is the same as in~\eqref{eqn:SemismoothNewtonStateEquation}, and as described in~\cite{Reyes2011} the derivative of $\regNonDifferentiability^\regParamMax(\strain(\velocity))$ is given as
\begin{align*}
	\begin{aligned}
		d^G_{(\velocity,\pressure)} & \left( \regNonDifferentiability^\regParamMax(\strain(\velocity)) \right) [(\perturbvelocity, \perturbpressure)] \\
		=& \frac{\yieldStress \regParam \, \strain(\perturbvelocity) : \strain(\adjvelocity)}{\max\nolimits_\regParamMax(\yieldStress, \regParam\! \left\| \strain(\velocity) \right\|)} \\
		&- \begin{cases}
			\bm{0} &\text{if } \yieldStress - \frac{1}{2 \regParamMax} \geq \regParam \! \left\| \strain(\velocity) \right\|, \\
			\frac{\yieldStress \regParam^2 \regParamMax \left( \strain(\velocity) : \strain(\perturbvelocity) \right) \, \strain(\velocity)}{ \max\nolimits_\regParamMax(\yieldStress, \regParam\! \left\| \strain(\velocity) \right\|)^2 \left\| \strain(\velocity) \right\|} \left( \regParam \left\| \strain(\velocity) \right\| - g + \frac{1}{2 \regParamMax} \right) & \text{if } \left| \regParam \! \left\| \strain(\velocity) \right\| - \yieldStress \right| \leq \frac{1}{2 \regParamMax}, \\
			\frac{\yieldStress (\strain(\velocity) : \strain(\perturbvelocity)) \, \strain(\velocity)}{\left\| \strain(\velocity) \right\|^3} &  \text{if } \yieldStress + \frac{1}{2 \regParamMax} \leq \regParam \! \left\| \strain(\velocity) \right\|.
		\end{cases}
	\end{aligned}
\end{align*}
Similar to \eqref{eqn:NonlinearTermStateEquation} we define a functional that describes the second term
\begin{align*}
	&\bm{d}_{\regParamMax} (\strain(\velocity), \strain(\perturbvelocity)) \\
	&= \begin{cases}
			\bm{0} &\text{if } \yieldStress - \frac{1}{2 \regParamMax} \geq \regParam \! \left\| \strain(\velocity) \right\|, \\
			\frac{\yieldStress \regParam^2 \regParamMax \left( \strain(\velocity) : \strain(\perturbvelocity) \right) \, \strain(\velocity)}{ \max\nolimits_\regParamMax(\yieldStress, \regParam\! \left\| \strain(\velocity) \right\|)^2 \left\| \strain(\velocity) \right\|} \left( \regParam \left\| \strain(\velocity) \right\| - g + \frac{1}{2 \regParamMax} \right) & \text{if } \left| \regParam \! \left\| \strain(\velocity) \right\| - \yieldStress \right| \leq \frac{1}{2 \regParamMax}, \\
			\frac{\yieldStress (\strain(\velocity) : \strain(\perturbvelocity)) \, \strain(\velocity)}{\left\| \strain(\velocity) \right\|^3} &  \text{if } \yieldStress + \frac{1}{2 \regParamMax} \leq \regParam \! \left\| \strain(\velocity) \right\|.
	\end{cases}
\end{align*}
The Newton's method that solves for $(\perturbvelocity^{k+1},\perturbpressure^{k+1})$ reads
\begin{align}
	\label{eqn:SemismoothNewtonStateEquationRegularized}
	\begin{aligned}
		&\int_{\compDomain} 2 \viscosity \, \strain(\perturbvelocity^{k+1}) : \strain(\adjvelocity) + \density \, ((\perturbvelocity^{k+1} \cdot \nabla) \velocity^k) \cdot \adjvelocity + \density \, ((\velocity^k \cdot \nabla) \perturbvelocity^{k+1}) \cdot \adjvelocity   \\
		&\hphantom{\int_{\compDomain} \ } - \perturbpressure^{k+1} \Div(\adjvelocity) + \adjpressure \Div(\perturbvelocity^{k+1}) + \frac{\yieldStress \regParam \, \strain(\perturbvelocity^{k+1}) : \strain(\adjvelocity)}{\max(\yieldStress, \regParam\! \left\| \strain(\velocity^k) \right\|)} \\
		&\hphantom{\int_{\compDomain} \ } - \bm{d}_{\regParamMax} (\strain(\velocity^k), \strain(\perturbvelocity^{k+1})) : \strain(\adjvelocity) \dx \\
		&= \int_{\compDomain} 2 \viscosity \, \strain(\velocity^k) : \strain(\adjvelocity) + \density \, ((\velocity^k \cdot \nabla) \velocity^k) \cdot \adjvelocity - \pressure^k \Div(\adjvelocity) + \adjpressure \Div(\velocity^k) %
		- \bm{f} \cdot \adjvelocity \\
		&\hphantom{= \int_{\compDomain} \ }+ \frac{\yieldStress \regParam \, \strain(\velocity^k) : \strain(\adjvelocity)}{\max_\regParamMax(\yieldStress, \regParam\! \left\| \strain(\velocity^k) \right\|)} \dx  \quad \forall \adjvelocity, \adjpressure \in \mathcal{H}_0 \times L^2(\compDomain).
	\end{aligned}
\end{align}
Due to the regularization of the $\max$-operator, this is a linear equation in $(\perturbvelocity^{k+1},\perturbpressure^{k+1})$. Using the Newton update
$$\left( \velocity^{k+1}, \pressure^{k+1} \right) = \left( \velocity^k, \pressure^k \right) + \alpha^{k+1} \left( \perturbvelocity^{k+1},\perturbpressure^{k+1} \right)$$
with Armijo backtracking line search to set $\alpha^{k+1}$, the regularized state equation \eqref{eqn:OptimizationProblemPDEConstraintRegularized} can be solved.

\emph{Regularized adjoint equation.}
In a similar manner as before in the unregularized case, however with the Fréchet operator with respect to $(\velocity, \pressure)$ of $\regNonDifferentiability^\regParamMax$ applied to the direction $(\perturbvelocity, \perturbpressure)$ instead of the Gâteaux semiderivative of $\regNonDifferentiability$, we obtain the regularized adjoint equation 
\begin{align}
	\label{eqn:GateauxAdjointEquationRegularized}
	\begin{aligned}
		0
		&= \int_{\compDomain} \viscosity \, \nabla \perturbvelocity : \nabla \velocity \dx + \int_{\compDomain} 2 \viscosity \, \strain(\perturbvelocity) : \strain(\adjvelocity) + \density \, ((\perturbvelocity \cdot \nabla) \velocity) \cdot \adjvelocity + \density \, ((\velocity \cdot \nabla) \perturbvelocity) \cdot \adjvelocity  \\
		&\hphantom{=\int_{\compDomain} \viscosity \, \nabla \perturbvelocity : \nabla \velocity \dx + \int_{\compDomain}\ } - \perturbpressure \Div(\adjvelocity) + \adjpressure \Div(\perturbvelocity) + \frac{\yieldStress \regParam \, \strain(\perturbvelocity) : \strain(\adjvelocity)}{\max\nolimits_\regParamMax(\yieldStress, \regParam\! \left\| \strain(\velocity) \right\|)} \\
		&\hphantom{=\int_{\compDomain} \viscosity \, \nabla \perturbvelocity : \nabla \velocity \dx + \int_{\compDomain}\ } - \bm{d}_{\regParamMax} (\strain(\velocity), \strain(\perturbvelocity)) : \strain(\adjvelocity) \dx
	\end{aligned}
\end{align}
Solving this equation for $\adjvelocity, \adjpressure \in \mathcal{H}_0 \times  L^2(\compDomain)$ for all $ \perturbvelocity, \perturbpressure \in \mathcal{H}_0 \times L^2(\compDomain)$ yields the regularized adjoint.

\emph{Regularized shape derivative.}
Again, by using the regularized $\max$-operator, we obtain the regularized shape derivative as the Eulerian derivative, where only the Eulerian derivative of the now-regularized term changes from \eqref{eqn:NonlinearTermDesignEquationM} to
\begin{align*}
	&\bm{m}_{\regParamMax} (\velocity, \strain(\dot{\velocity}), \meshDeformation) \\
	&= \begin{cases}
		\bm{0} & \text{if } \yieldStress - \frac{1}{2 \regParamMax} \geq \regParam \! \left\| \strain(\velocity) \right\|, \\
		\begin{aligned}[c]
		&\yieldStress \tfrac{\regParam^2 \regParamMax \, \strain(\velocity) \, \left( \strain(\velocity) : \left( \strain(\dot{\velocity}) - \tfrac{1}{2} \left( \nabla \velocity \nabla \meshDeformation + {\nabla \meshDeformation}^\top {\nabla \velocity}^\top \right) \right) \right)}{\max\nolimits_{\regParamMax}(\yieldStress, \regParam\! \left\| \strain(\velocity) \right\|)^2 \left\| \strain(\velocity) \right\|} \\
		&\qquad \cdot \left( \regParam \left\| \strain(\velocity) \right\| - g + \tfrac{1}{2 \regParamMax} \right)
		\end{aligned} 
		& \text{if } \left| \regParam \! \left\| \strain(\velocity) \right\| - \yieldStress \right| \leq \frac{1}{2 \regParamMax}, \\
		\yieldStress \frac{\strain(\velocity)}{\left\| \strain(\velocity) \right\|^3} \, \strain(\velocity) : \left( \strain(\dot{\velocity}) - \tfrac{1}{2} \left( \nabla \velocity \nabla \meshDeformation + {\nabla \meshDeformation}^\top {\nabla \velocity}^\top \right) \right) & \text{if } \yieldStress + \frac{1}{2 \regParamMax} \leq \regParam \! \left\| \strain(\velocity) \right\|.
		\end{cases}
\end{align*}
After including other differentiable terms and inserting the weak form of the regularized state~\eqref{eqn:OptimizationProblemPDEConstraintRegularized} and of the adjoint~\eqref{eqn:GateauxAdjointEquationRegularized} this yields
\begin{align}
	\label{eqn:RegularizedShapeDerivative}
	\begin{aligned}
		&d_E L_A^\regParamMax(\shape, (\velocity,\pressure),(\adjvelocity,\adjpressure)) [\meshDeformation] \\
		&= \int_{\compDomain} - \viscosity \left( \nabla \velocity \nabla \meshDeformation \right) : \nabla \velocity
		- \viscosity \left( \nabla \velocity \nabla \meshDeformation + {\nabla \meshDeformation}^\top {\nabla \velocity}^\top  \right) : \strain(\adjvelocity) \\
		&\hphantom{=\int_{\compDomain} \,} - \viscosity \, \strain(\velocity) : \left( \nabla \adjvelocity \nabla \meshDeformation + {\nabla \meshDeformation}^\top {\nabla \adjvelocity}^\top  \right) + \bm{q}_{\regParamMax} (\velocity, \meshDeformation) : \strain(\adjvelocity) \\ \\
		&\hphantom{=\int_{\compDomain} \,} - \frac{\yieldStress \regParam \left( \nabla \velocity \nabla \meshDeformation + {\nabla \meshDeformation}^\top {\nabla \velocity}^\top  \right) : \strain(\adjvelocity)}{2\max\nolimits_{\regParamMax}(\yieldStress, \regParam\! \left\| \strain(\velocity) \right\|)} - \frac{\yieldStress \regParam \, \strain(\velocity) : \left( \nabla \adjvelocity \nabla \meshDeformation + {\nabla \meshDeformation}^\top {\nabla \adjvelocity}^\top  \right)}{2\max\nolimits_{\regParamMax}(\yieldStress, \regParam\! \left\| \strain(\velocity) \right\|)} \\
		&\hphantom{=\int_{\compDomain} \,} - \density \left( \nabla \meshDeformation \velocity \cdot \nabla \right) \velocity \cdot \adjvelocity + \pressure {\nabla \adjvelocity}^\top : \nabla \meshDeformation - \left(\nabla \bm{f} \, \meshDeformation \right) \cdot \adjvelocity - \adjpressure {\nabla \velocity}^\top : \nabla \meshDeformation \\
		&\hphantom{=\int_{\compDomain} \,} + \Div(\meshDeformation) \left( \frac{\viscosity}{2} \nabla \velocity : \nabla \velocity + 2 \viscosity \, \strain(\velocity) : \strain(\adjvelocity)  + \frac{\yieldStress \regParam \, \strain(\velocity) : \strain(\adjvelocity)}{\max\nolimits_{\regParamMax}(\yieldStress, \regParam\! \left\| \strain(\velocity) \right\|)} \right. \\
		&\hphantom{=\int_{\compDomain} + \Div(\meshDeformation) \ \ }  \left. \vphantom{\frac{\strain(\adjvelocity)}{\max(\left\| \strain(\velocity) \right\|)}} + \density \left(\velocity \cdot \nabla \right) \velocity \cdot \adjvelocity - \pressure \Div(\adjvelocity) - \bm{f} \cdot \adjvelocity + \adjpressure \Div(\velocity) \right) \dx \\
		&\hphantom{= \ } + \left( \lambda_1 - \nu \left(\vol(\shape) - \mathcal{V} \right) \right) \int_{\shape} \meshDeformation^\top \bm{n} \,\mathrm{d}s \\
		&\hphantom{= \ } + \frac{1}{\vol(\shape)} \left( \nu \left(\bary(\shape) - \bm{\mathcal{B}})^\top  - \left(\lambda_2, \lambda_3 \right) \right) \right) \int_{\shape} (\bary(\shape)-\bm{x}) \, \meshDeformation^\top \bm{n} \,\mathrm{d}s \\
		&\hphantom{= \ } + \left( \lambda_4 - \nu \left(\peri(\shape) - \mathcal{P} \right) \right) \int_{\shape} \bm{n}^\top {\nabla \meshDeformation} \bm{n} - \Div(\meshDeformation) \,\mathrm{d}s,
	\end{aligned}
\end{align}
where we have the regularized version of \eqref{eqn:NonlinearTermDesignEquationQ} given as
\begin{align*}
	\begin{aligned}
		&\bm{q}_{\regParamMax} (\velocity, \meshDeformation)
		&= \begin{cases}
			\bm{0} & \text{if } \yieldStress - \frac{1}{2 \regParamMax} \geq \regParam \! \left\| \strain(\velocity) \right\|, \\
			\begin{aligned}[c]
			&\yieldStress \tfrac{\regParam^2 \regParamMax \, \strain(\velocity) \, \left( \strain(\velocity) : \left( \nabla \velocity \nabla \meshDeformation + {\nabla \meshDeformation}^\top {\nabla \velocity}^\top \right) \right)}{2 \max\nolimits_{\regParamMax}(\yieldStress, \regParam\! \left\| \strain(\velocity) \right\|)^2 \left\| \strain(\velocity) \right\|} \\
			&\qquad \cdot \left( \regParam \left\| \strain(\velocity) \right\| - g + \tfrac{1}{2 \regParamMax} \right)
			\end{aligned} 
			& \text{if } \left| \regParam \! \left\| \strain(\velocity) \right\| - \yieldStress \right| \leq \frac{1}{2 \regParamMax}, \\
			\yieldStress \frac{\strain(\velocity)}{2 \left\| \strain(\velocity) \right\|^3} \, \strain(\velocity) : \left( \nabla \velocity \nabla \meshDeformation + {\nabla \meshDeformation}^\top {\nabla \velocity}^\top \right) & \text{if } \yieldStress + \frac{1}{2 \regParamMax} \leq \regParam \! \left\| \strain(\velocity) \right\|.
		\end{cases}
	\end{aligned}
\end{align*}
With the state, adjoint and design equation in unregularized as well as regularized form, we are now ready to describe the algorithm.

\subsection{Optimization algorithm}
\label{sec:OptimizationAlgorithm}

In order to optimize on the shape space $M_1^\text{ad}$ we need
\begin{itemize}
	\item  a Riemannian metric to define the search direction %
	with respect to the corresponding metric and 
	\item a method of mapping from the tangent space of the manifold to the manifold itself in order to formulate the next (shape-)iterate in an algorithm. 
\end{itemize} 

The  exponential map is the first option for a mapping from the tangent space of the manifold to the manifold itself, however due to the numerical effort of solving a second-order differential equation for the exponential map, it is common to use a retraction instead. 
In our experiments in \Cref{sec:Num}, we concentrate on the retraction already used in \cite{GeiersbachSuchanWelker_NS,Pryymak2023,Suchan2023}.
This retraction updates the shape in each iteration and
 is closely related to the perturbation of  identity, which 
 is defined for vector fields on the domain.
 Let $\compDomain$ be the considered domain. If we consider $\compDomain$ depending on a shape $\shape\subset\compDomain$, the perturbation of identity acting on $\compDomain$ in the direction of a sufficiently smooth vector field $\bm{V} \colon \compDomain \to \mathbb{R}^2$ gives
\begin{equation}
\label{UpdateDomain}
\compDomain_t(u)= \{ \bm{x} \in \compDomain \colon \bm{x} + t \bm{V} (\bm{x}) \}
\end{equation}
with %
$t\geq 0$.
Here, the corresponding  shape update for  $\shape \subset \compDomain$ is given by
$\shape_t = \{ \bm{x} \in \shape \colon \bm{x} + t\bm{v}(\bm{x}) \},$
where $\bm{v}=\bm{V}|_{\shape}$.

As mentioned in \Cref{sub:ShapeSpace}, we concentrate on the Steklov-Poincaré metric firstly defined in \cite{Schulz2016a}. %
If we consider this metric and define $\Xi \coloneqq \{ \bm{W} \in H^1(\compDomain, \mathbb{R}^2) \colon \bm{W} = \bm{0} \text{ on } \partial \compDomain \setminus \shape\}$, the above-mentioned vector field  $\meshDeformation$ %
is obtained by solving the deformation equation
\begin{align}
	\label{eqn:deformationEquation}
	a(\meshDeformation, \bm{W}) = -\widetilde{d_E L_A}(\shape, (\velocity, \pressure), (\adjvelocity, \adjpressure))[\bm{W}] \quad \forall \bm{W} \in \Xi
\end{align}
in the unregularized case and 
\begin{align}
	\label{eqn:deformationEquationRegularized}
	a(\meshDeformation, \bm{W}) = -d_E L_A^\regParamMax(\shape, (\velocity, \pressure), (\adjvelocity, \adjpressure))[\bm{W}] \quad \forall \bm{W} \in \Xi
\end{align}
in the regularized case, where $a$ is a symmetric and coercive bilinear form in both cases.

As in \cite{LoayzaGeiersbachWelkerHandbook,GeiersbachSuchanWelker_NS,geiersbach2023stochastic,Pryymak2023} the optimization algorithm on $\admissibleShapeSpace$ can be formulated as in \Cref{alg:algorithmOptManifold}, in which we directly formulate the outcome of using the above-mentioned retraction related to the perturbation of identity in line 13.
For the regularized $\max$-operator, this algorithm can be used directly. Without regularization of the $\max$-operator, we additionally implement a safeguard using the regularized problem (similar to~\cite[Algorithm~2]{Luft2020}) in case no step size that satisfies the Armijo condition can be found.

\begin{algorithm}
	\caption{Augmented Lagrangian method on $\admissibleShapeSpace$}
	\label{alg:algorithmOptManifold}
	\begin{algorithmic}[1]
		\State \textbf{Input:} Initial shape $\shape \in \admissibleShapeSpace$ with $\shape\subset \compDomain$, maximum step size $t_{\max}>0$, $\tau \in (0,1)$, $\xi>1$ %
		\State \textbf{Initialization:} $\bm{\lambda} \in \mathbb{R}^4$, $\nu > 0$
		\While{$\bm{\lambda}$ not converged}
		\State $c_1 = \left\| \bm{c}(\shape) \right\|$
		\While{$\shape$ not converged}
		\State Solve state equation~\eqref{eqn:OptimizationProblemPDEConstraintOptimalityCondition} or~\eqref{eqn:OptimizationProblemPDEConstraintRegularized} to obtain $(\velocity, \pressure) \in \mathcal{H} \times L^2(\compDomain)$
		\State Compute objective functional value~$L_A$ in~\eqref{eqn:augmentedLagrangeFunctional} or $L_A^\regParamMax$ in~\eqref{eqn:augmentedLagrangeFunctionalRegularized}
		\State Solve adjoint equation~\eqref{eqn:GateauxAdjointEquationApprox} or~\eqref{eqn:GateauxAdjointEquationRegularized} to obtain $(\adjvelocity , \adjpressure) \in \mathcal{H}_0 \times L^2(\compDomain)$
		\State \parbox[t]{298pt}{Determine potential descent direction~$\meshDeformation$ by solving~\ref{eqn:deformationEquation} with right-hand side~\ref{eqn:GateauxShapeDerivativeApprox} or by solving~\ref{eqn:deformationEquationRegularized} with right-hand side~\eqref{eqn:RegularizedShapeDerivative}}
		\State Determine step size $t$ by Armijo backtracking with $t_{\max}$ on $L_A$ or $L_A^\regParamMax$
		\If{No valid step size~$t$ found}
		\State \parbox[t]{283pt}{Determine state, adjoint, shape gradient and step size~$t$ with regularized model}%
		\EndIf
		\State %
		Update $\shape$ by updating $\compDomain$ (cf. \eqref{UpdateDomain})
		\EndWhile
		\State $\bm{\lambda} \leftarrow \bm{\lambda} - \nu \, \bm{c}(\shape)$
		\If{$\left\| \bm{c}(\shape) \right\| \geq \tau c_1$}
		\State $\nu \leftarrow \xi \nu$
		\EndIf
		\EndWhile
	\end{algorithmic}
\end{algorithm}

\section{Numerical results}
\label{sec:Num}

The numerical investigation is performed using dolfinx version 0.6.0 \cite{Scroggs2022,Scroggs2022a,Alnaes2014} on a regular workstation. The workstation has 32 cores and 256 gigabytes of RAM, however no parallelization is used. The mesh of the unit domain $\compDomain=(0,1)^2$ (the interior of the shape $\shape$, which is a square an edge length of $0.1$ centered at $(0.3, 0.45)^\top$, is not part of $\compDomain$) is generated using Gmsh \cite{Geuzaine2009} with a mesh size of $h=\frac{1}{40}$ with local refinement around the square to $h=\frac{1}{200}$, which yields 4963 nodes and 9606 triangular elements. 
The bilinear form in \eqref{eqn:deformationEquation} is chosen as linear elasticity, cf., e.g., \cite{Siebenborn2017,Geiersbach2021,GeiersbachSuchanWelker_NS}, i.e.,
\begin{align*}
	a(\meshDeformation, \bm{W}) = \int_\compDomain 2 \hat{\mu} \, \strain(\meshDeformation) : \strain(\bm{W}) \dx,
\end{align*}
where $\hat{\mu}$ is obtained by interpolating $\hat{\mu}=5$ on $\shape$ and $\hat{\mu}=1$ on $\partial \compDomain \setminus \shape$ into the computational domain $\compDomain$ by solving a Poisson problem with Dirichlet boundary conditions.

To ensure inf-sup stability, Taylor-Hood $P2-P1$ elements are used for the discretization of the velocity and the pressure, cf., e.g., \cite{Taylor1974,Bercovier1979}, and $P1$ elements are used for $\meshDeformation$. The initial penalty factor for the geometrical constraints is chosen as $\nu=10^5$, the Lagrange multipliers are initialized as $\bm{\lambda}=\bm{0}$, the fluid viscosity is $\viscosity=1$, volumetric forces were neglected $\bm{f}=\bm{0}$, the yield threshold is $\yieldStress=20$, and the regularization parameter is set to $\regParam=10^3$. The regularization of the $\max$-operator is done using $\regParamMax=10^{-1}$. The augmented Lagrange parameters are chosen as $\tau=0.9$ and $\xi=2$. We set the geometrical equality constraints to $\mathcal{V}=0.04$, $\bm{\mathcal{B}}=\left(0.3, 0.45\right)^\top$ and $\mathcal{P}=0.76$. The velocity on the left boundary is given by a parabola with a value of $1$ at the center and $0$ at the corners of the domain, i.e., $\velocity=(-4 x_2(x_2-1),0)^\top$, and zero on all other parts of $\Gamma^D$. The outflow is modeled by a homogenous Neumann condition on the right-hand boundary of the domain. The maximum step size for the backtracking is given as $t_{\max} = 6.25 \cdot 10^{-6}$. Each inner loop (cf. \Cref{alg:algorithmOptManifold}, line~5--13) is executed for $2\,000$ iterations, while the outer loop (line~3--16) is stopped after $10$ iterations. Therefore, we obtain $20\,000$ cumulative iterations in total.

\subsection{Numerical solution of the state equation and active set}
\label{sec:NumericalSolutionStateEquationAndActiveSet}

\begin{figure}[tb]
	\centering
	\setlength\figureheight{7cm} 
	\setlength\figurewidth{.48\textwidth}
	\includegraphics{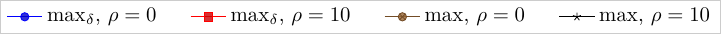}\\%
	\includetikz{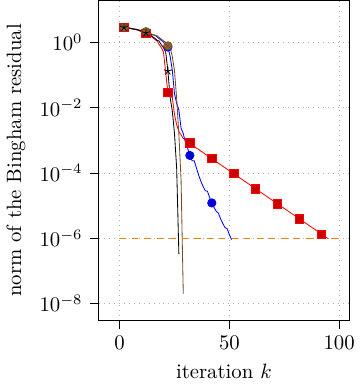}%
	\includetikz{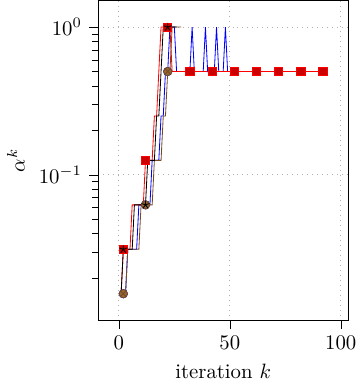}%
	\caption{Residual norm and accepted step sizes of the iterative scheme using $\rho=0$ and $\rho=10$ with regularized and unregularized $\max$-operator for the solution of the state equation.}
	\label{fig:NumericalResultsStateSolver}
\end{figure}

In a first step, we compare the solvers of the state equation in the regularized and unregularized case. No good initial guess for the solution of the state equation is available, therefore we consider an initial guess obtained by solving the Stokes equation: Find $\velocity, \pressure \in \mathcal{H} \times L^2(\compDomain)$ s.t.
\begin{align*}
	\label{eqn:OptimizationProblemPDEConstraintInitialGuess}
	\begin{aligned}
		0 = \int_{\compDomain} &2 \viscosity \, \strain(\velocity) : \strain(\adjvelocity) - \pressure \Div(\adjvelocity) + \adjpressure \Div(\velocity)- \bm{f} \cdot \adjvelocity \dx \quad \forall \adjvelocity, \adjpressure \in \mathcal{H}_0 \times L^2(\compDomain).
	\end{aligned}
\end{align*}
Plots of the norm of the residual of the discretized Bingham equation over iterations~$k$ and of the accepted step sizes $\alpha^k$ of the iterative schemes to obtain the state after an initial guess from the Stokes equation can be found in \Cref{fig:NumericalResultsStateSolver}. We observe faster than linear convergence with the unregularized $\max$-operator, which is basically unaffected by the presence of the convective term. Full steps are taken in the last iterations. The regularized $\max$-operator on the other hand does not always use full Newton steps. Without the convective term ($\rho=0$) we observe a quicker reduction in residual for the regularized $\max$ operator than with the convective term ($\rho=10$) due to the stronger nonlinearity in $(\velocity,\pressure)$.

\begin{figure}[tb]
	\centering
	\setlength\figureheight{.45\textwidth} 
	\setlength\figurewidth{.45\textwidth}
	\includegraphics[width=\textwidth]{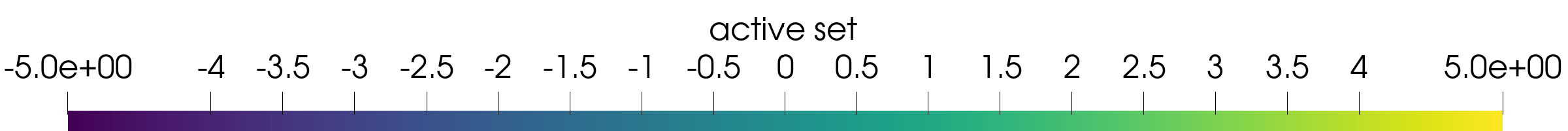}\\%
	\includegraphics[width=\figurewidth]{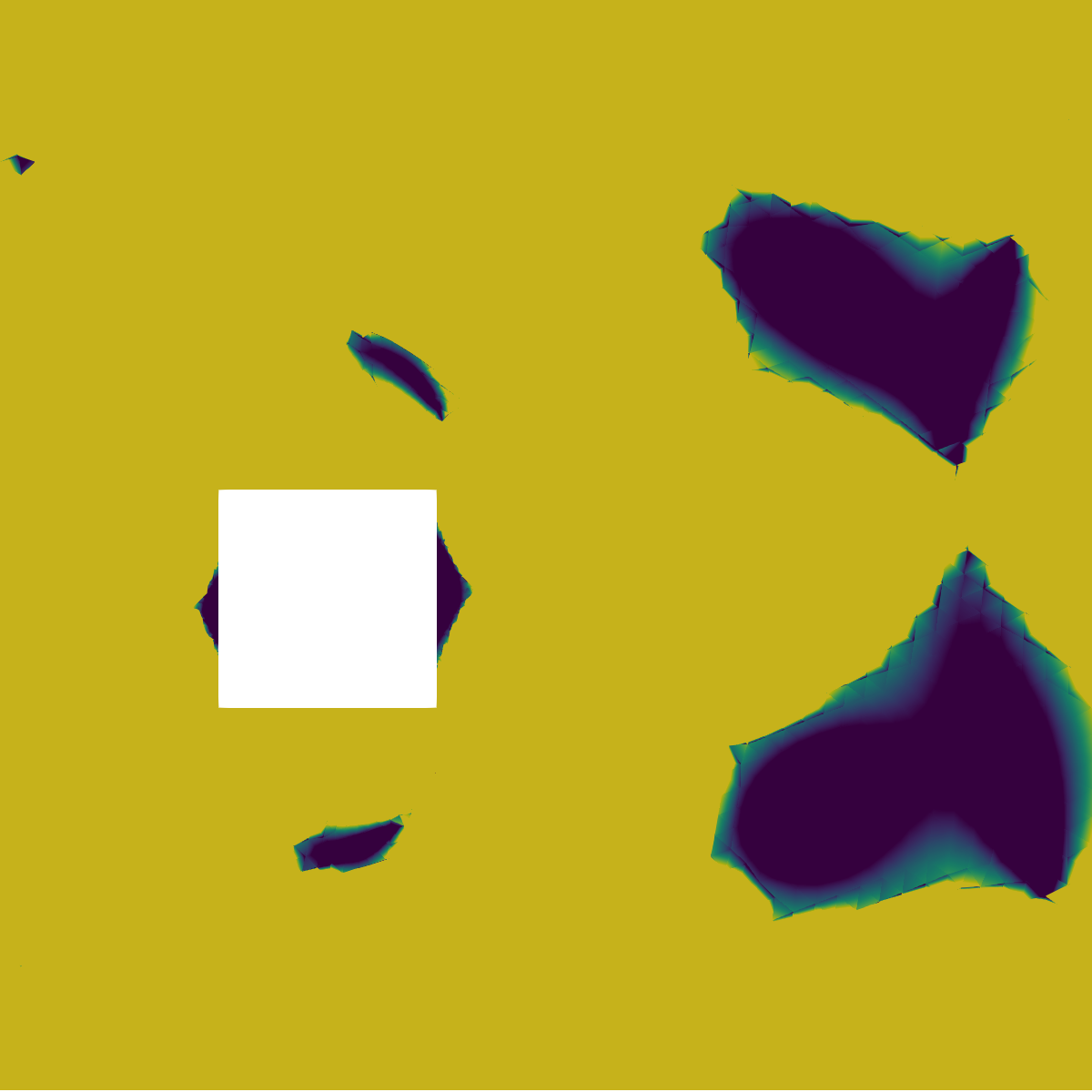}%
	\quad%
	\includegraphics[width=\figurewidth]{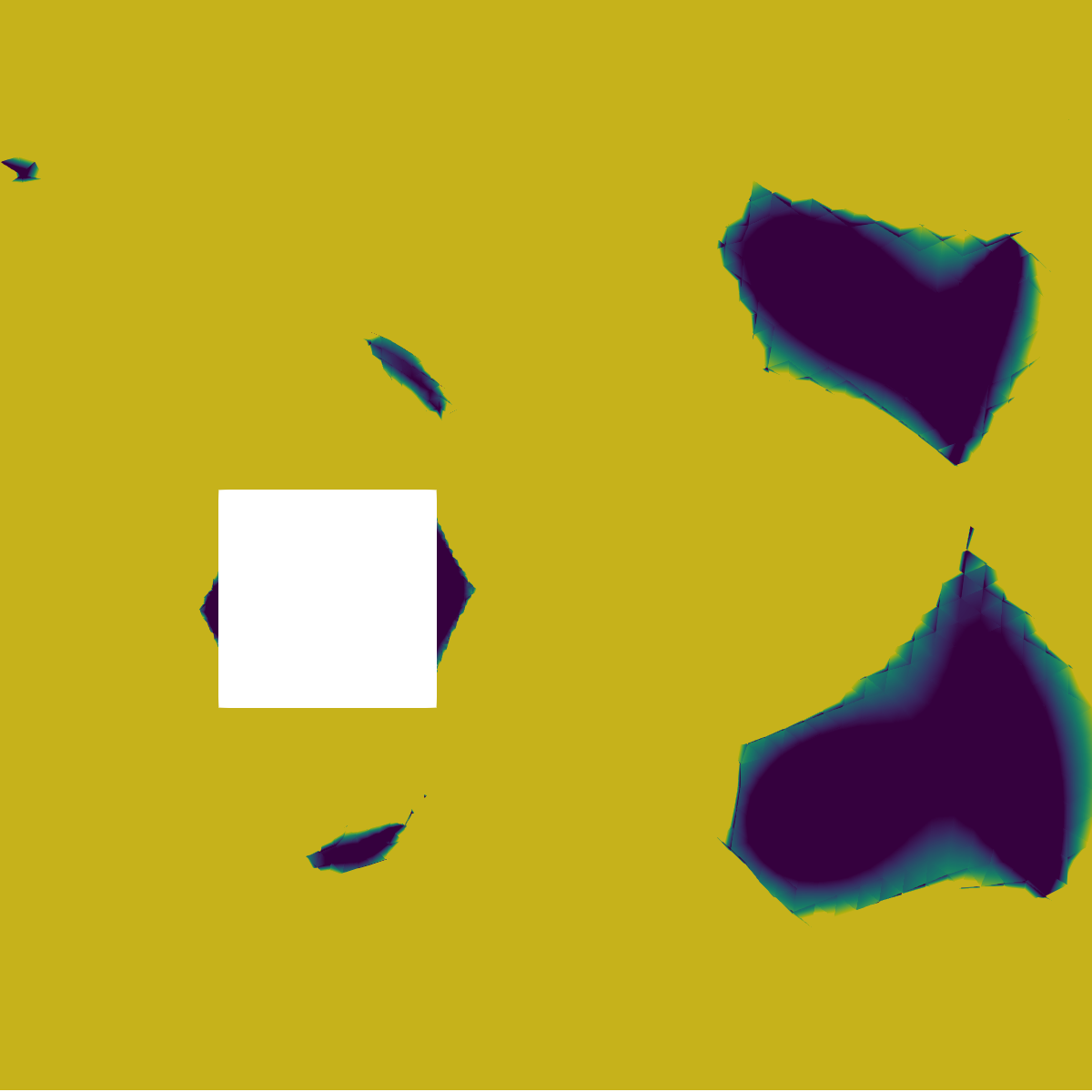}%
	\caption{Active set $\regParam \! \left\| \strain(\velocity) \right\| - \yieldStress$ at the start of the optimization with regularization of the $\max$-operator and $\rho=0$ (left) and $\rho=10$ (right).}
	\label{fig:NumericResults_activeSet_regularized}
\end{figure}

\begin{figure}[tb]
	\centering
	\setlength\figureheight{.45\textwidth} 
	\setlength\figurewidth{.45\textwidth}
	\includegraphics[width=\textwidth]{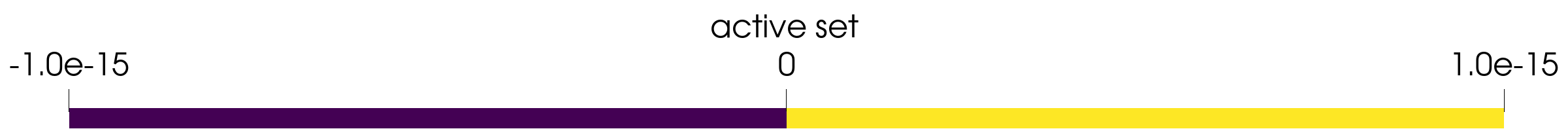}\\%
	\includegraphics[width=\figurewidth]{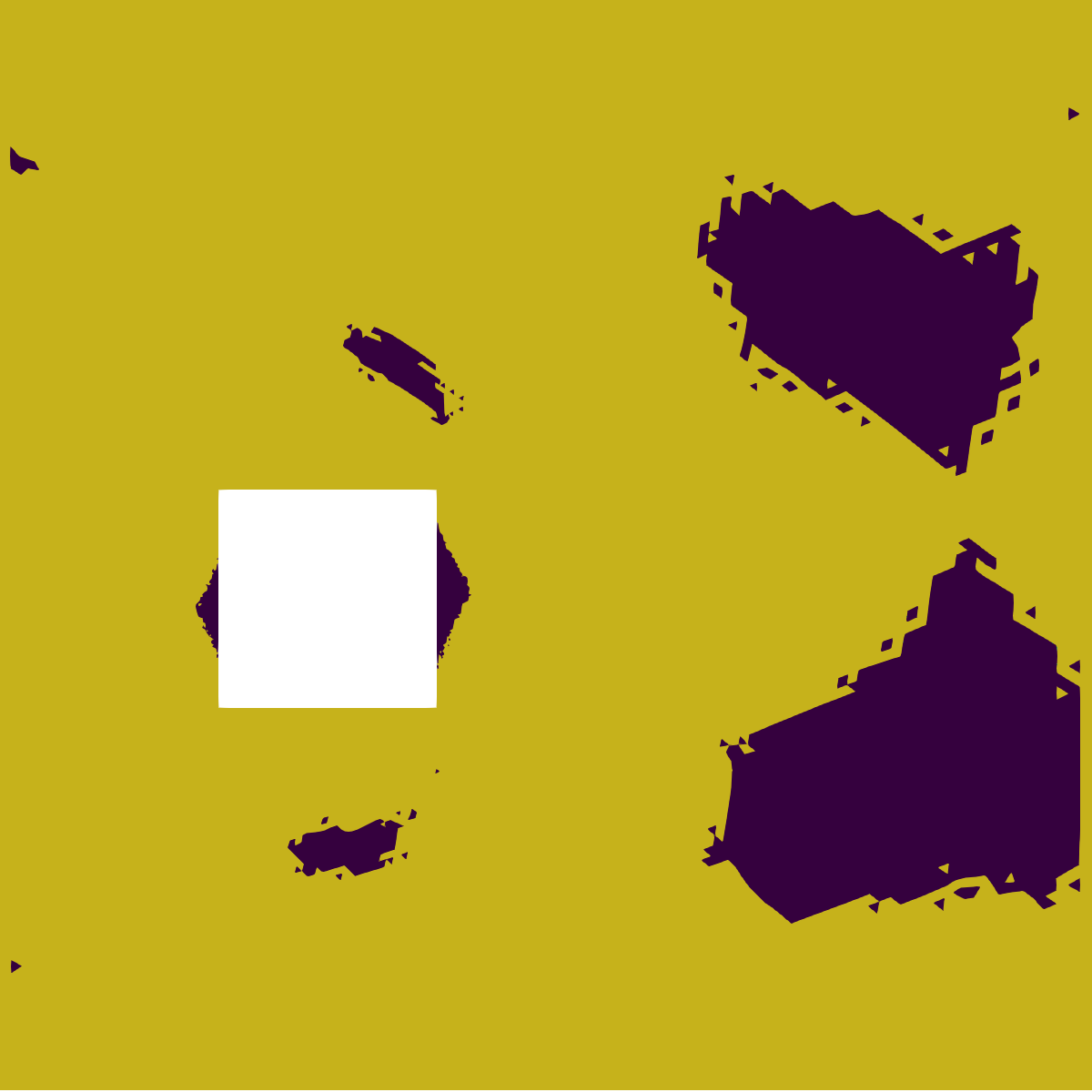}%
	\quad%
	\includegraphics[width=\figurewidth]{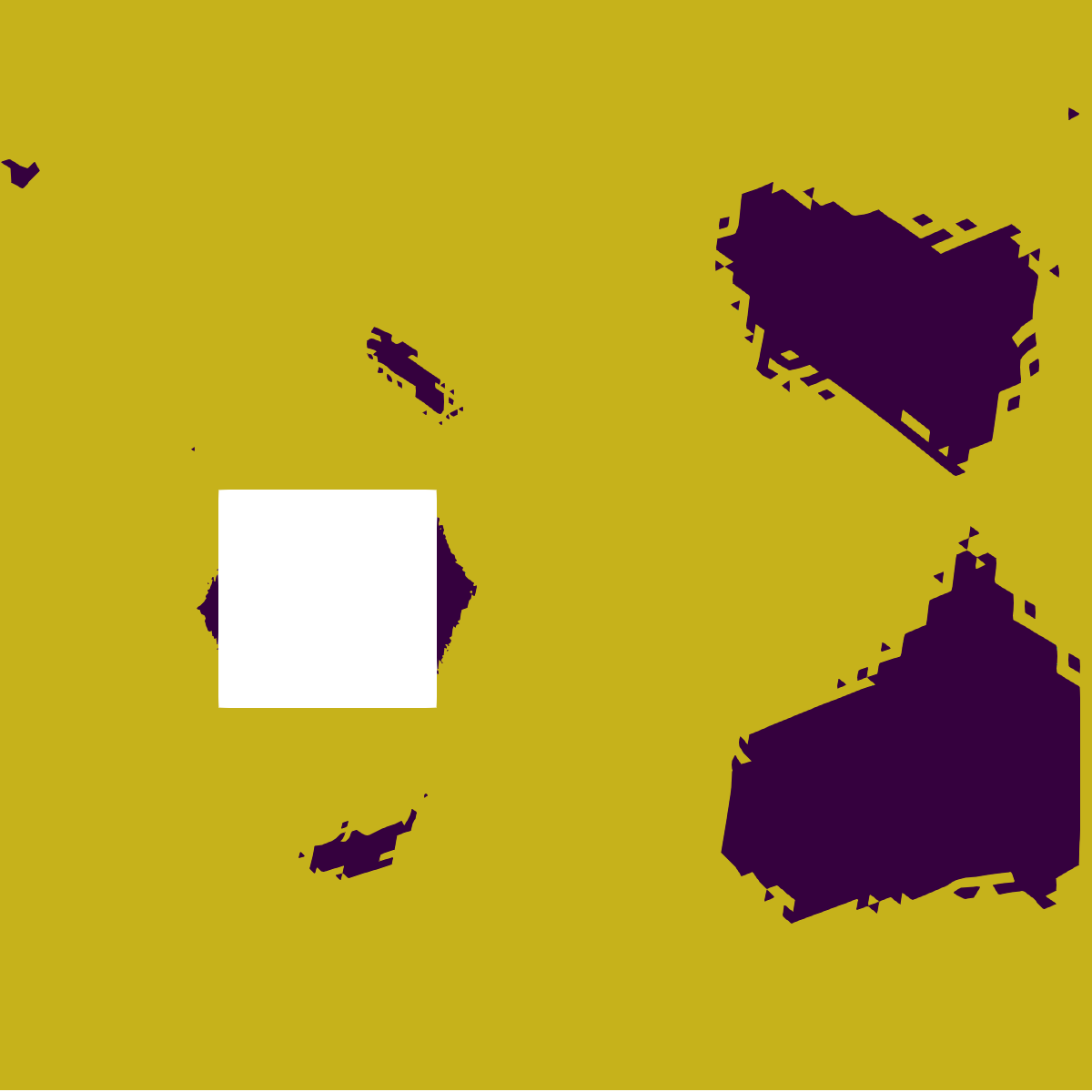}%
	\caption{Active set $\regParam \! \left\| \strain(\velocity) \right\| - \yieldStress$ at the start of the optimization without regularization of the $\max$-operator and $\rho=0$ (left) and $\rho=10$ (right).}
	\label{fig:NumericResults_activeSet_unregularized}
\end{figure}

In \Cref{fig:NumericResults_activeSet_regularized} we provide an illustration of the parts of the domain with regularization in which $\yieldStress - \frac{1}{2 \regParamMax} \geq \regParam \! \left\| \strain(\velocity) \right\|$ (blue), $\left| \regParam \! \left\| \strain(\velocity) \right\| - \yieldStress \right| \leq \frac{1}{2 \regParamMax}$ (intermediate colors) and $\yieldStress + \frac{1}{2 \regParamMax} \leq \regParam \! \left\| \strain(\velocity) \right\|$ (yellow). These illustrations are obtained by projecting the functional $\regParam \! \left\| \strain(\velocity) \right\| - \yieldStress$ to a discontinuous linear finite element space. Values which exceed the colorbar are colored as the color for the minimum and maximum value, respectively. We observe small differences between the domains.
In \Cref{fig:NumericResults_activeSet_unregularized}, a similar illustration is given, however now in the unregularized case we only distinguish the two cases $\yieldStress > \regParam \! \left\| \strain(\velocity) \right\|$ (blue) and $\yieldStress \leq \regParam \! \left\| \strain(\velocity) \right\|$ (yellow). Similar differences as in the unregularized case are observed with respect to the convective term. The domain is similarly split into the different sets. We observe small areas around larger parts of the domain where $\yieldStress > \regParam \! \left\| \strain(\velocity) \right\|$, which we attribute to the errors introduced by the projection to a finite element space.

\subsection{Shape optimization without convective term ($\rho=0$)}

\begin{figure}[tbp]
	\centering
	\setlength\figureheight{7cm} 
	\setlength\figurewidth{.48\textwidth}
	\includetikz{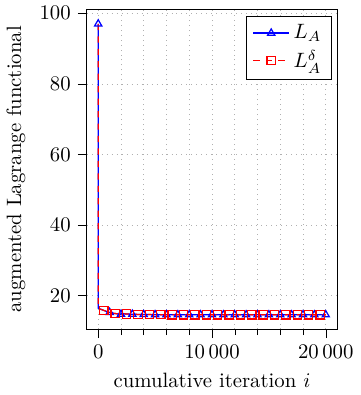}%
	\includetikz{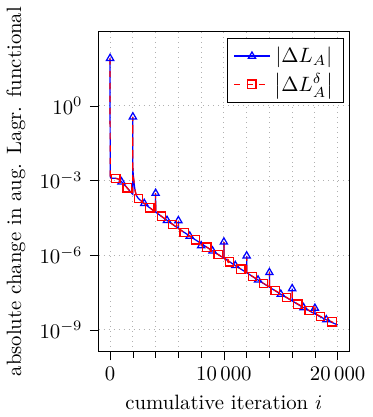}%
	\caption{Augmented Lagrange functional (left) and absolute change of augmented Lagrange functional $\Delta (\cdot ) = \left|(\cdot )^{i}-(\cdot )^{i-1}\right|$, $i=1,\ldots,20\,000$, (right) using $\rho=0$ with regularized and unregularized $\max$-operator.}
\label{fig:NumericResults_L_A_Stokes}
\end{figure}

\begin{figure}[tbp]
	\centering
	\setlength\figureheight{7cm} 
	\setlength\figurewidth{\textwidth}
	\includetikz{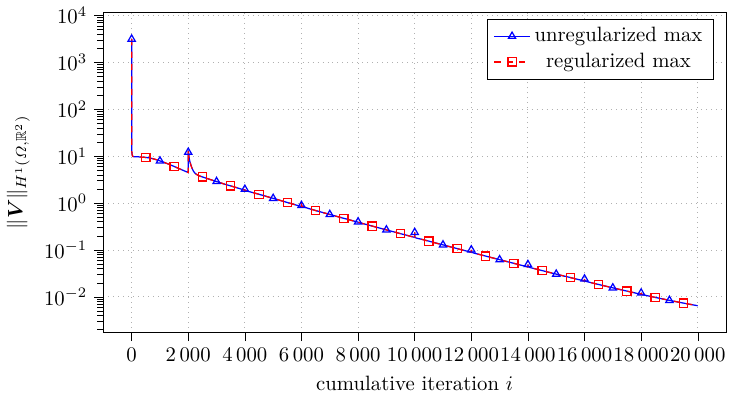}%
	\caption{$H^1$-norm of the mesh deformation using $\rho=0$ with regularized and unregularized $\max$-operator.}
	\label{fig:NumericResults_normV_Stokes}
\end{figure}

\begin{figure}[tbp]
	\centering
	\setlength\figureheight{.45\textwidth} 
	\setlength\figurewidth{.45\textwidth}
	\includegraphics[width=\textwidth]{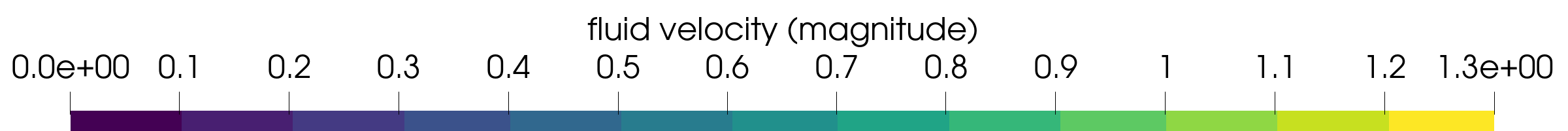}\\%
	\includegraphics[width=\figurewidth]{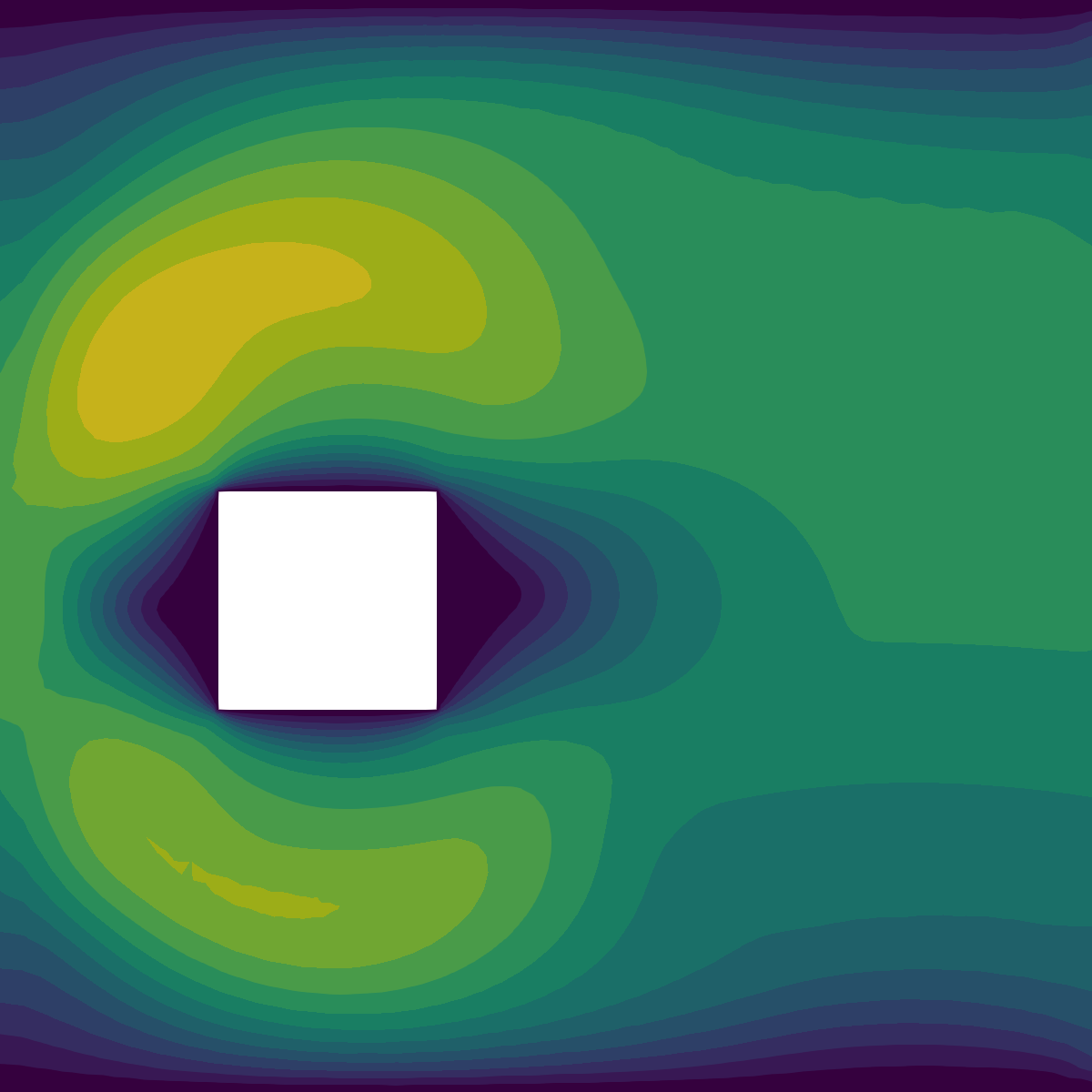}%
	\quad%
	\includegraphics[width=\figurewidth]{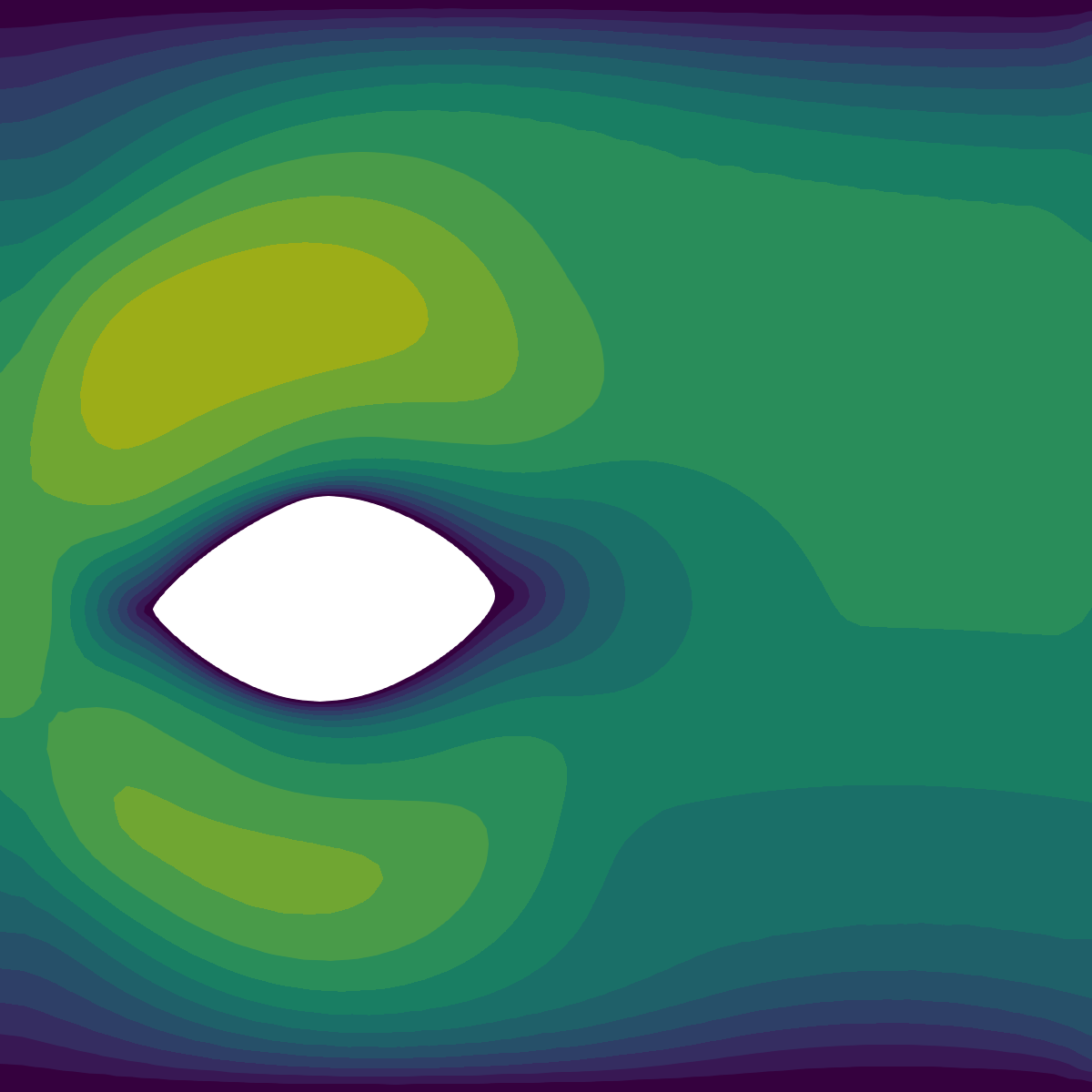}%
	\caption{Fluid velocity magnitude at the start of the optimization (left) and at the end of the optimization (right) with regularization of the $\max$-operator and $\rho=0$.}
	\label{fig:NumericResults_shapes_Stokes_regularized}
\end{figure}

\begin{figure}[tbp]
	\centering
	\setlength\figureheight{.45\textwidth} 
	\setlength\figurewidth{.45\textwidth}
	\includegraphics[width=\textwidth]{figs/no_convective_term/colorbar}\\%
	\includegraphics[width=\figurewidth]{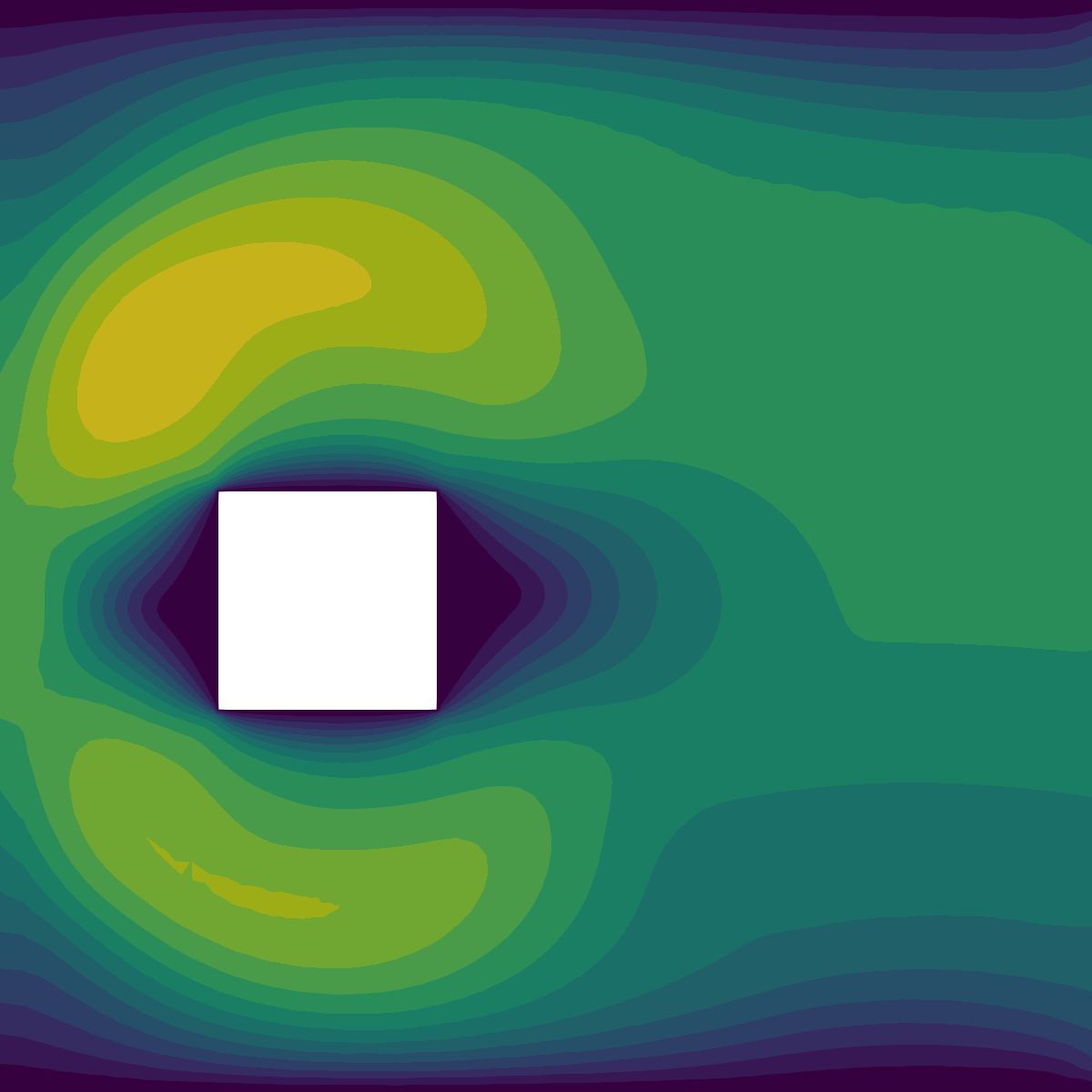}%
	\quad%
	\includegraphics[width=\figurewidth]{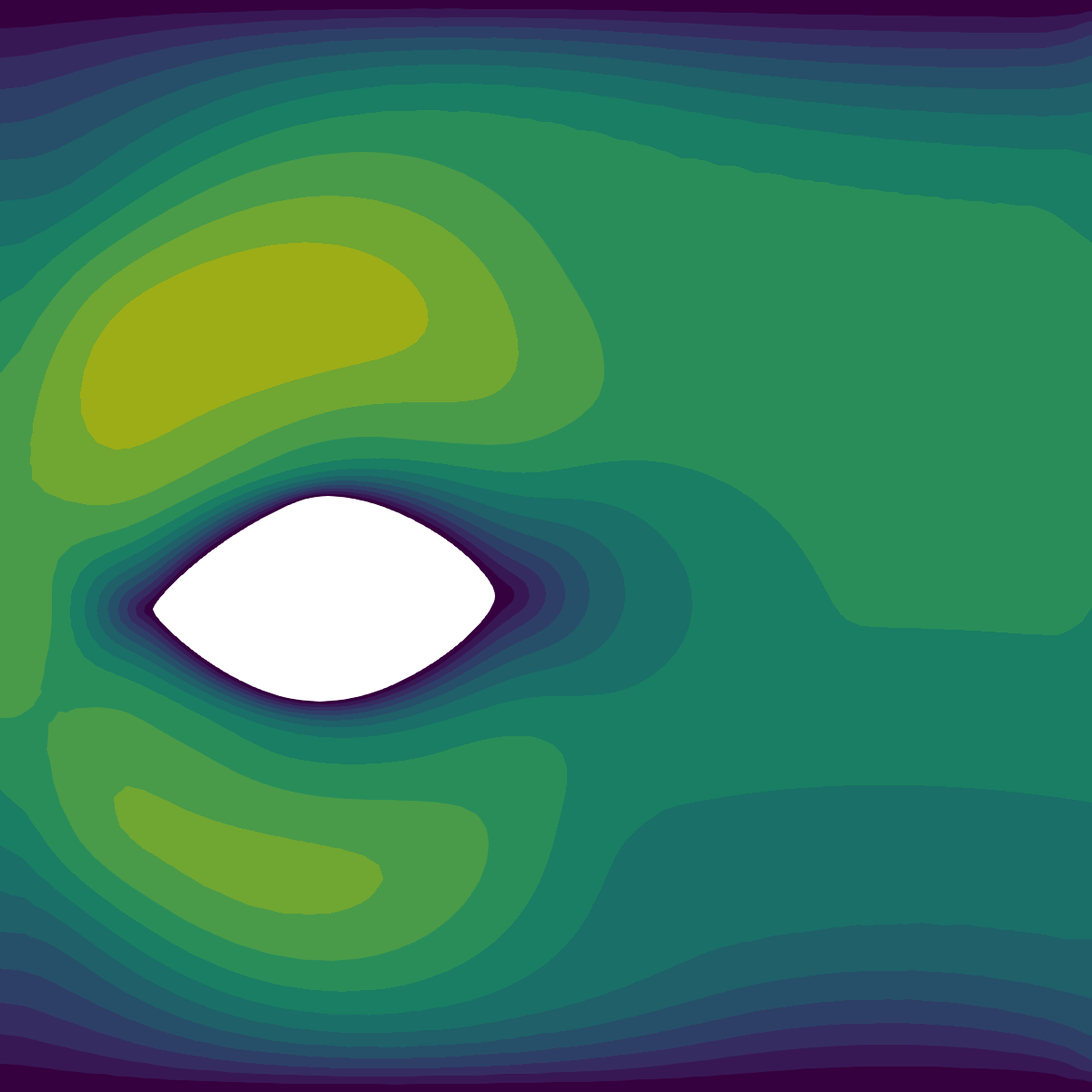}%
	\caption{Fluid velocity magnitude at the start of the optimization (left) and at the end of the optimization (right) without regularization of the $\max$-operator and $\rho=0$.}
	\label{fig:NumericResults_shapes_Stokes_unregularized}
\end{figure}

We now apply the optimization algorithm from \Cref{sec:OptimizationAlgorithm}.
For the first optimization we neglect the influence of the convective term, similar to \cite{Reyes2011}, i.e., $\rho=0$, and consider the problem with the regularized as well as the unregularized $\max$-operator. To accelerate the computation of the state after the initial computation as in \Cref{sec:NumericalSolutionStateEquationAndActiveSet}, the state of the previous shape optimization iteration is used as an initial guess. The state with a regularized $\max$-operator converges rather slowly in the range from $10^{-4}$ to $10^{-6}$ in comparison to an unregularized $\max$-operator (cf. \Cref{fig:NumericalResultsStateSolver}). For the second optimization iteration the iterative solver of the state with $\rho=0$ required only $17$ iterations until a residual below $10^{-6}$ is reached again in the unregularized case. In the regularized case, the Newton solver required $38$ iterations. The plots of the augmented Lagrange functional and change in augmented Lagrange functional between two subsequent optimization iterations can be found in \Cref{fig:NumericResults_L_A_Stokes} for both cases. Furthermore, the $H^1$-norm of the mesh deformation is shown in~\Cref{fig:NumericResults_normV_Stokes}, and the shapes at the start and at the end of the optimization are shown in \Cref{fig:NumericResults_shapes_Stokes_regularized} in the regularized and in \Cref{fig:NumericResults_shapes_Stokes_unregularized} in the unregularized case. We observe a strong initial reduction in augmented Lagrange functional value, mainly caused by the initial violation of the geometrical constraints. The stops of the inner gradient descent loop are clearly visible every $2\,000$ iterations in both, the change in augmented Lagrange functional as well as mesh deformation. A reduction by more than five orders of magnitude is obtained for $\| \meshDeformation \|$. No difference is visible in any of the figures nor in the obtained shapes between the optimization with a regularized $\max$-operator and an unregularized $\max$-operator.

\subsection{Shape optimization with convective term ($\rho=10$)}

\begin{figure}[tbp]
	\centering
	\setlength\figureheight{7cm} 
	\setlength\figurewidth{.48\textwidth}
	\includetikz{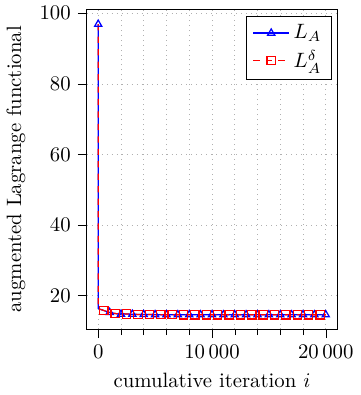}%
	\includetikz{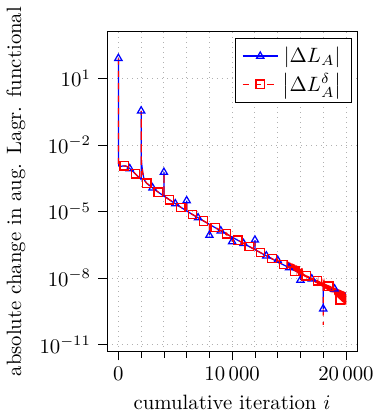}%
	\caption{Augmented Lagrange functional (left) and absolute change of augmented Lagrange functional $\Delta (\cdot ) = \left|(\cdot )^{i}-(\cdot )^{i-1}\right|$, $i=1,\ldots,20\,000$, (right) using $\rho=10$ with regularized and unregularized $\max$-operator.}
	\label{fig:NumericResults_L_A_NavierStokes}
\end{figure}

\begin{figure}[tbp]
	\centering
	\setlength\figureheight{7cm} 
	\setlength\figurewidth{\textwidth}
	\includetikz{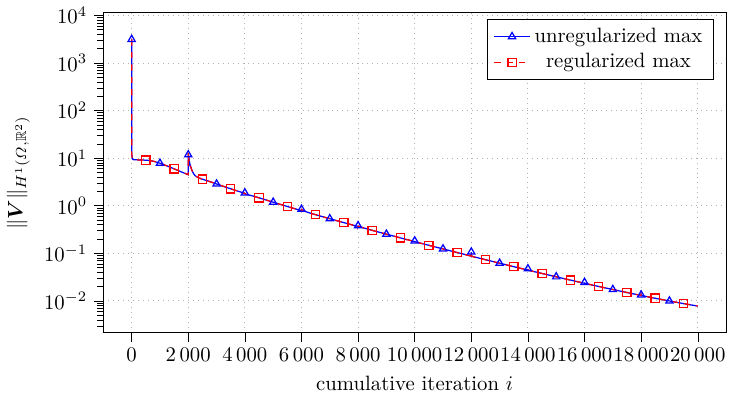}%
	\caption{$H^1$-norm of the mesh deformation using $\rho=10$ with regularized and unregularized $\max$-operator.}
	\label{fig:NumericResults_normV_NavierStokes}
\end{figure}

\begin{figure}[tbp]
	\centering
	\setlength\figureheight{.45\textwidth} 
	\setlength\figurewidth{.45\textwidth}
	\includegraphics[width=\textwidth]{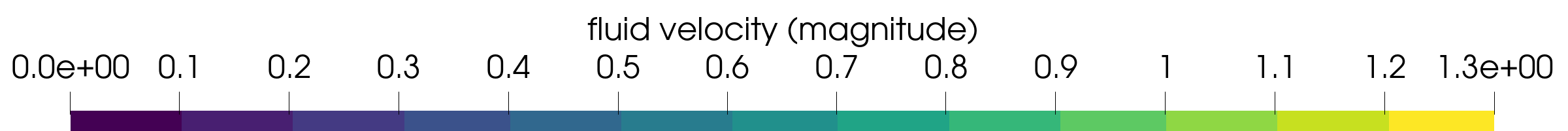}\\%
	\includegraphics[width=\figurewidth]{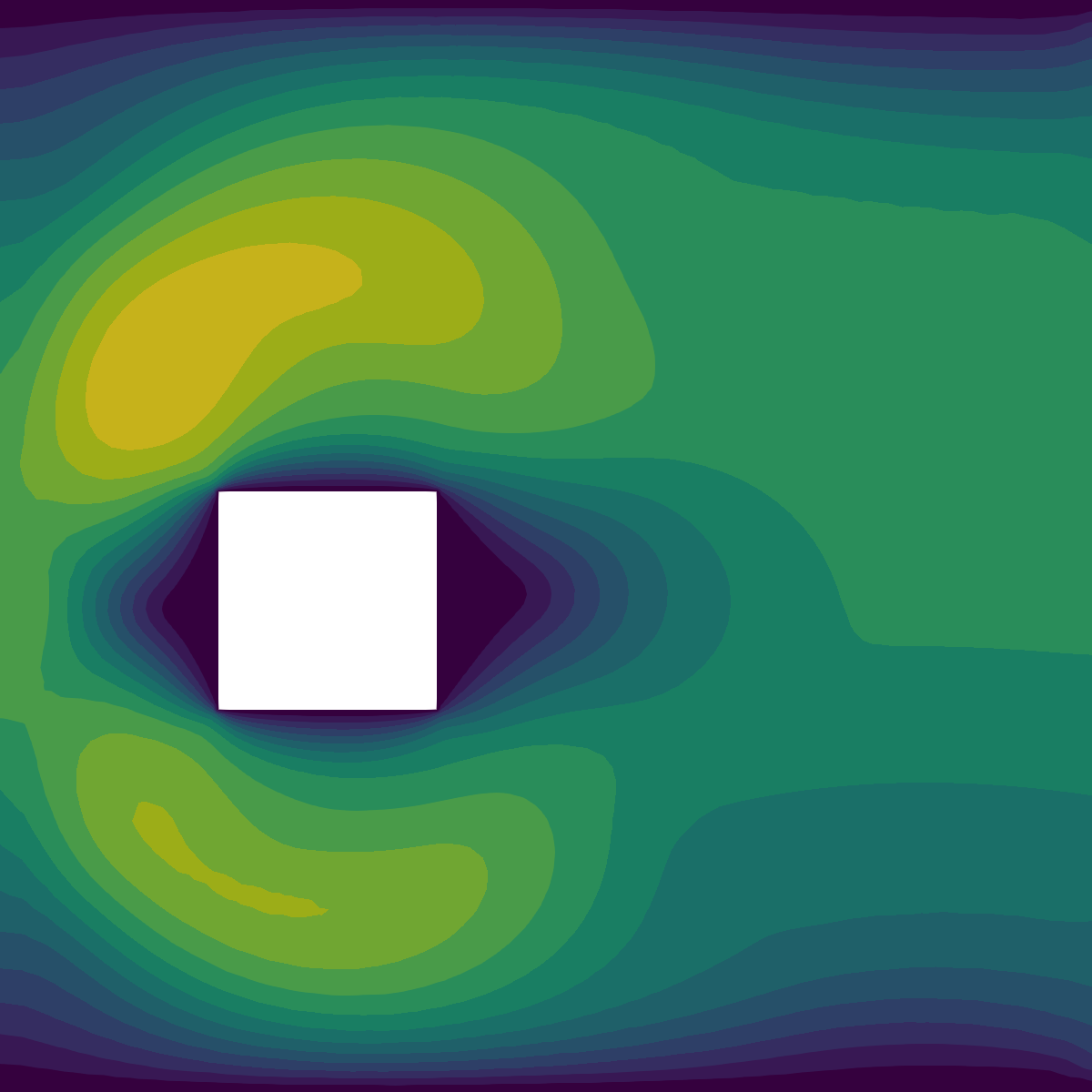}%
	\quad%
	\includegraphics[width=\figurewidth]{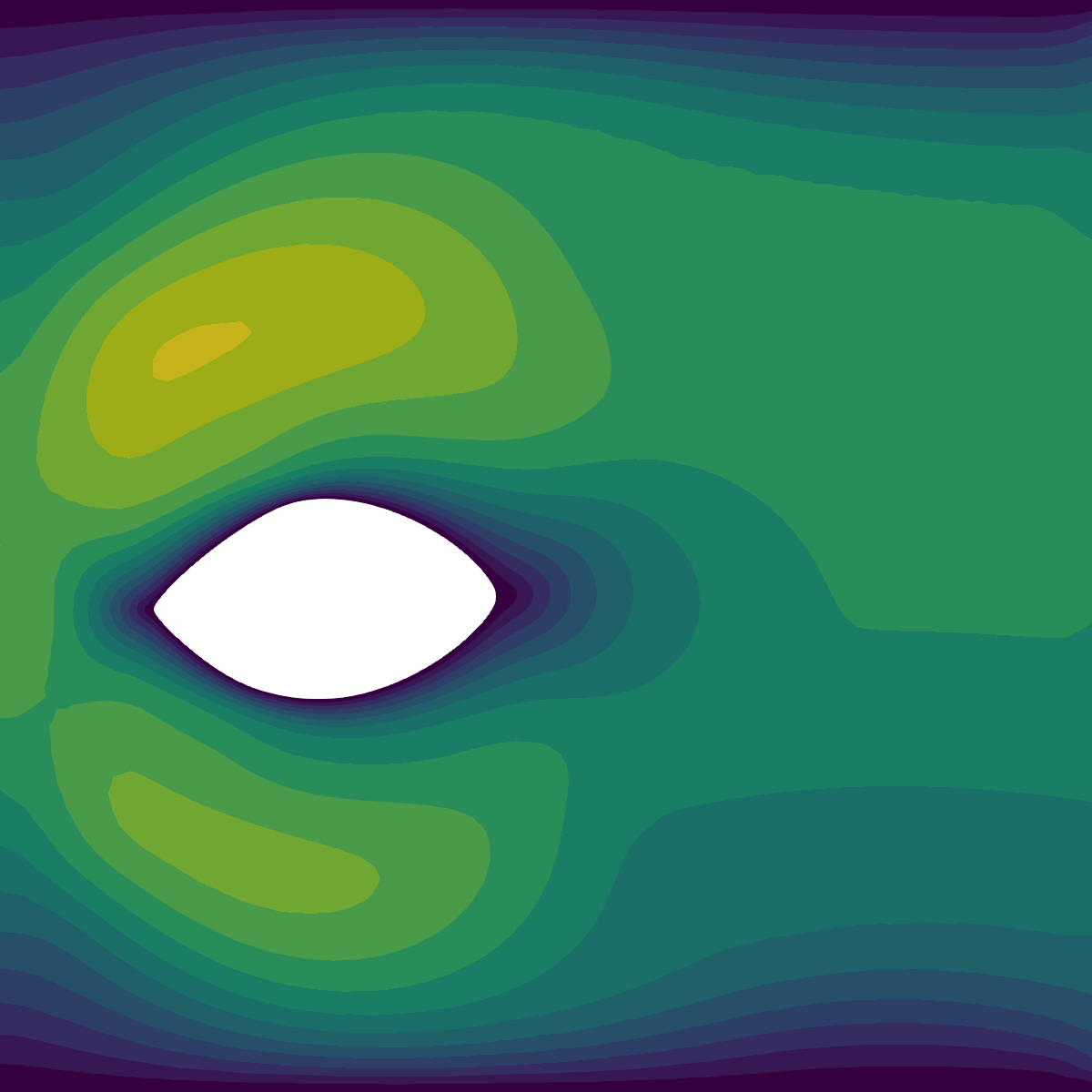}%
	\caption{Fluid velocity magnitude at the start of the optimization (left) and at the end of the optimization (right) with regularization of the $\max$-operator and $\rho=10$.}
	\label{fig:NumericResults_shapes_NavierStokes_regularized}
\end{figure}

\begin{figure}[tbp]
	\centering
	\setlength\figureheight{.45\textwidth} 
	\setlength\figurewidth{.45\textwidth}
	\includegraphics[width=\textwidth]{figs/convective_term/colorbar}\\%
	\includegraphics[width=\figurewidth]{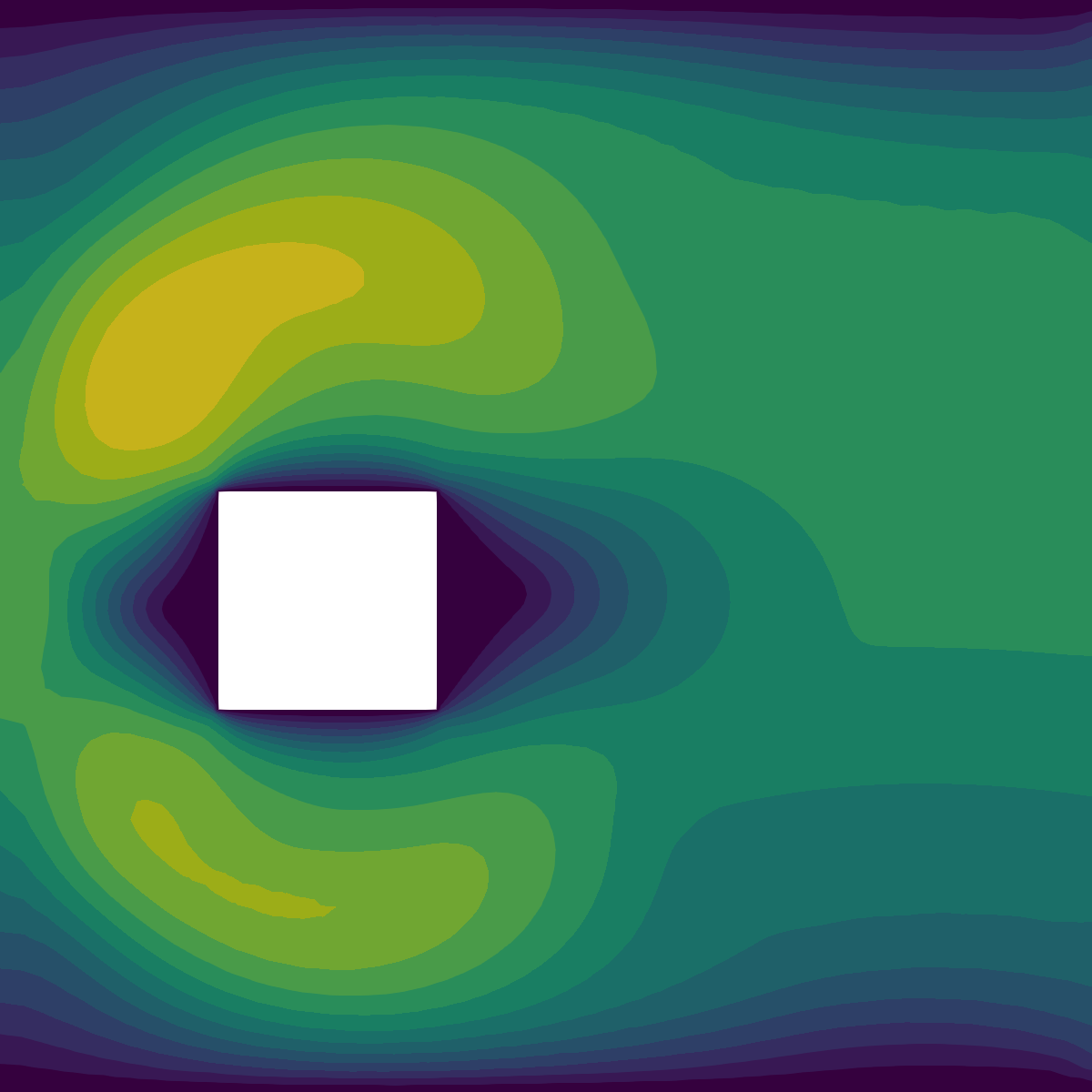}%
	\quad%
	\includegraphics[width=\figurewidth]{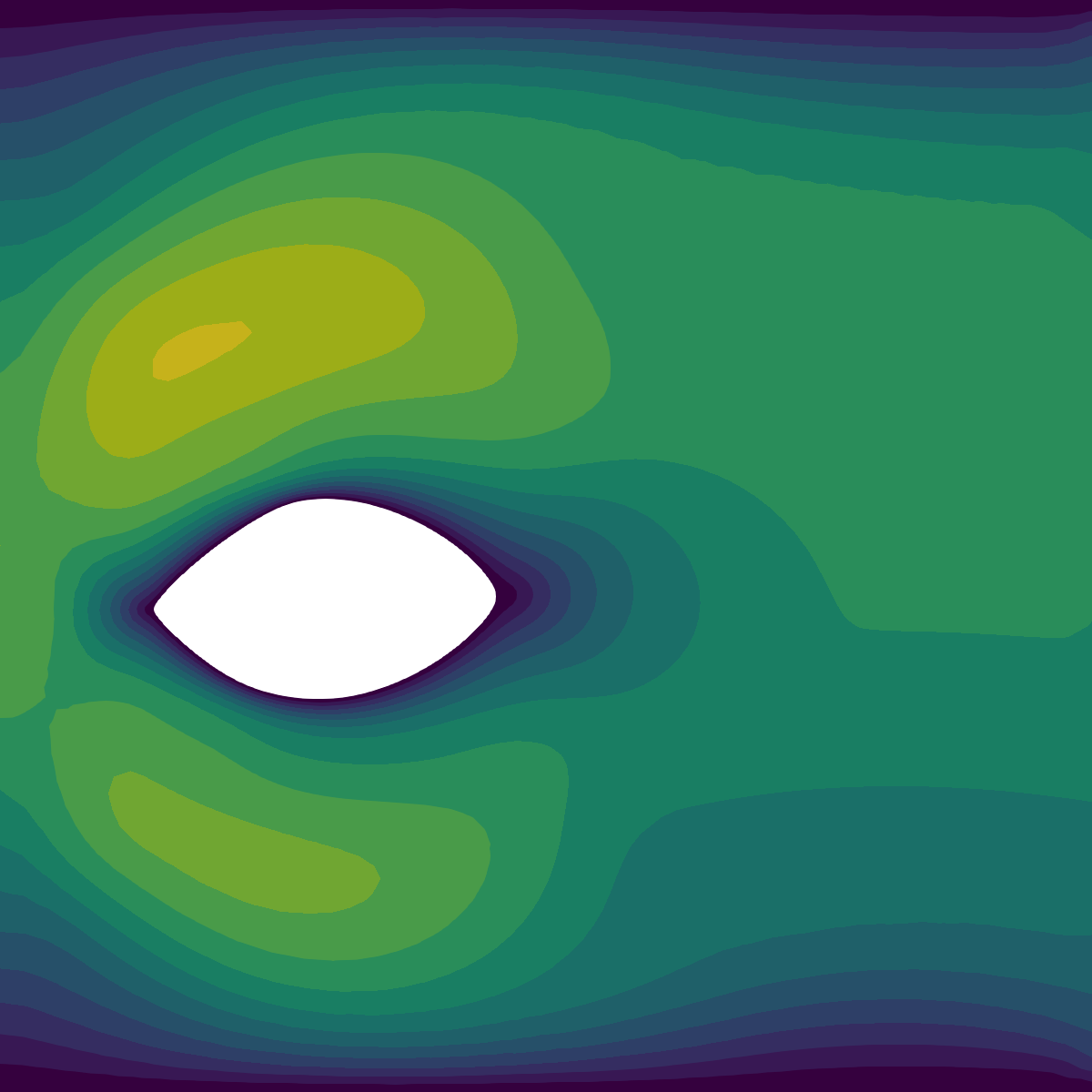}%
	\caption{Fluid velocity magnitude at the start of the optimization (left) and at the end of the optimization (right) without regularization of the $\max$-operator and $\rho=10$.}
	\label{fig:NumericResults_shapes_NavierStokes_unregularized}
\end{figure}

A second experiment is performed that includes the influence of the convective term. Here, we use $\rho = 10$. As before, the plots of the augmented Lagrange functional and change in augmented Lagrange functional between two subsequent iterations are given in \Cref{fig:NumericResults_L_A_NavierStokes} for both, the regularized and the unregularized case. The $H^1$-norm of the mesh deformation is provided in~\Cref{fig:NumericResults_normV_NavierStokes}, and the shapes before and after optimization are shown in \Cref{fig:NumericResults_shapes_NavierStokes_regularized} in the regularized and in \Cref{fig:NumericResults_shapes_NavierStokes_unregularized} in the unregularized case. As in the example without the convective term we observe the same behavior. The results with and without regularization are nearly indistinguishable from each other. The colorbars are chosen the same as without the convective term. We observe a slightly different velocity field between $\rho=0$ and $\rho=10$. 

\section{Conclusion}
Bingham flow poses challenges not only to simulation but also to optimization, because of its variational inequality nature. In this paper, Gâteaux semiderivatives are used in order to derive a Newton-like solution strategy for the state equation as well necessary conditions for shape optimization based on adjoints. Because of the variational inequality aspect, the adjoint equation contains nonlinear terms, which are approximated linearly. The solution of the resulting approximate adjoint equation is used within a gradient based optimization algorithm. Because it is not guaranteed that the resulting step has a descent property, a safeguard technique based on a regularized problem formulation is employed, as in \cite{Luft2020}. We observe however that the safeguard has never become active and thus never had to save the step. Thus, further research seems necessary for the theoretical justification of the descent property of the approximate step.

\par\addvspace{\bigskipamount}

\begin{acknowledgements}
	\noindent \textbf{Acknowledgements} This work has been partly supported by the German Research Foundation (DFG) within the priority program SPP~1962 under contract number WE~6629/1-1 and SCHU~804/19-1. 
\end{acknowledgements}

\bibliographystyle{plainurl}
\bibliography{literature.bib}

\end{document}